\documentclass[twocolumn]{autart}    

\usepackage{amsmath, amssymb, mathtools, accents, dsfont}
\usepackage{algorithm,algpseudocode}
\usepackage{balance}
\usepackage{graphicx}
\usepackage{subfigure}
\usepackage{mathtools}
\usepackage{blkarray, bigstrut}
\usepackage{booktabs}
\usepackage{bbm}
\usepackage{tikz}
	\usetikzlibrary{shapes,arrows}
	\tikzstyle{frame} = [draw, -latex]
	\tikzstyle{line} = [draw]
	\tikzstyle{line2} = [draw, dashdotted]
	\tikzstyle{line3} = [draw, dashed]
	\tikzstyle{line3UD} = [draw, dashed]
	\tikzstyle{place} = [circle, draw=black, fill=white, thick, inner sep=2pt, minimum size=1mm]
	\tikzstyle{place2} = [circle, draw=black, fill=black, thick, inner sep=2pt, minimum size=1mm]
	\tikzstyle{placeRed} = [circle, draw=red, fill=red, thick, inner sep=2pt, minimum size=1mm]
	\tikzstyle{vertex} = [circle, draw=black, fill=black, thick, inner sep=2pt, minimum size=1mm]
\usepackage{tkz-euclide}

\usetikzlibrary{shapes,arrows} 

\tikzstyle{decision} = [diamond, draw, fill=blue!20,
    text width=4.5em, text badly centered, node distance=3cm, inner sep=0pt]
\tikzstyle{block1} = [rectangle, draw, text width=8em, text centered, minimum height=4em]
\tikzstyle{block2} = [rectangle, draw, text width=3em, text centered, minimum height=4em]
\tikzstyle{block3} = [rectangle, draw, text width=11em, text centered, minimum height=12em, dashed,black]
\tikzstyle{block4} = [rectangle, draw, text width=11em, text centered, minimum height=18em, dashed,black]
\tikzstyle{block5} = [rectangle, draw, text width=11em, text centered, minimum height=32em, dashed,black]
\tikzstyle{block6} = [rectangle, draw, text width=11em, text centered, minimum height=18.5em, dashed,black]
\tikzstyle{block7} = [rectangle, draw, text width=11em, text centered, minimum height=11.8em, dashed,black]
\tikzstyle{line01} = [draw, -latex']
\tikzstyle{line02} = [draw, latex'-latex']

\usepackage{setspace}
\usepackage[switch,pagewise]{lineno}
\newcommand*\patchAmsMathEnvironmentForLineno[1]{%
\expandafter\let\csname old#1\expandafter\endcsname\csname #1\endcsname
\expandafter\let\csname oldend#1\expandafter\endcsname\csname end#1\endcsname
\renewenvironment{#1}%
{\linenomath\csname old#1\endcsname}%
{\csname oldend#1\endcsname\endlinenomath}}%
\newcommand*\patchBothAmsMathEnvironmentsForLineno[1]{%
\patchAmsMathEnvironmentForLineno{#1}%
\patchAmsMathEnvironmentForLineno{#1*}}%
\AtBeginDocument{%
\patchBothAmsMathEnvironmentsForLineno{equation}%
\patchBothAmsMathEnvironmentsForLineno{align}%
\patchBothAmsMathEnvironmentsForLineno{flalign}%
\patchBothAmsMathEnvironmentsForLineno{alignat}%
\patchBothAmsMathEnvironmentsForLineno{gather}%
\patchBothAmsMathEnvironmentsForLineno{multline}%
}

\makeatletter
\def\BState{\State\hskip-\ALG@thistlm}
\makeatother
\algdef{SE}[DOWHILE]{Do}{doWhile}{\algorithmicdo}[1]{\algorithmicwhile\ #1}%

\algnewcommand\algorithmicswitch{\textit{switch}}
\algnewcommand\algorithmiccase{\textit{case}}
\algnewcommand\algorithmicassert{\texttt{assert}}
\algnewcommand\Assert[1]{\State \algorithmicassert(#1)}%
\algdef{SE}[SWITCH]{Switch}{EndSwitch}[1]{\algorithmicswitch\ #1\ }{\algorithmicend\ \algorithmicswitch}%
\algdef{SE}[CASE]{Case}{EndCase}[1]{\algorithmiccase\ #1}{\algorithmicend\ \algorithmiccase}%
\algtext*{EndSwitch}%
\algtext*{EndCase}%

\usepackage{amsmath,bm}

\begin{document}

\begin{frontmatter}
\title{Composition Rules for Strong Structural Controllability and Minimum Input Problem in Diffusively-Coupled Networks}


\author[Gwangju]{Nam-Jin Park}\ead{namjinpark@gist.ac.kr},    
\author[Gyeonggi]{Seong-Ho Kwon}\ead{seongho@krri.re.kr},               
\author[Daejeon1]{Yoo-Bin Bae}\ead{ybbae@kari.re.kr},               
\author[Daejeon2]{Byeong-Yeon Kim}\ead{byeongyeon@kaeri.re.kr},               
\author[Golden]{Kevin L. Moore}\ead{kmoore@mines.edu}, and
\author[Gwangju]{Hyo-Sung Ahn\corauthref{cor}}\ead{hyosung@gist.ac.kr}

\corauth[cor]{Corresponding author at: School of Mechanical Engineering, Gwangju Institute of Science and Technology (GIST), Gwangju, Korea}

\address[Gwangju]{School of Mechanical Engineering, Gwangju Institute of Science and Technology (GIST), Gwangju, Korea}                   
\address[Gyeonggi]{Korea Railroad Research Institute (KRRI), Gyeonggi-do, Korea}                   
\address[Daejeon1]{Korea Aerospace Research Institute (KARI), Daejeon, Korea}             
\address[Daejeon2]{Korea Atomic Energy Research Institute (KAERI), Daejeon, Korea}             
\address[Golden]{Department of Electrical Engineering, Colorado School of Mines, Golden, CO, USA}             


\begin{abstract}                  
This paper presents new results and reinterpretation of existing conditions for strong structural controllability in a structured network determined by the zero/non-zero patterns of edges. For diffusively-coupled networks with self-loops, we first establish a necessary and sufficient condition for strong structural controllability, based on the concepts of dedicated and sharing nodes. Subsequently, we define several conditions for strong structural controllability across various graph types by decomposing them into disjoint path graphs.  We further extend our findings by introducing a composition rule, facilitating the analysis of strong structural controllability in larger networks. This rule allows us to determine the strong structural controllability of connected graphs called pactus graphs (a generalization of the well-known cactus graph) by consideration of the strong structural controllability of its disjoint component graphs. In this process, we introduce the notion of a component input node, which is a state node that functions identically to an external input node.
Based on this concept, we present an algorithm with approximate polynomial complexity to 
determine the minimum number of external input nodes required to maintain strong structural controllability in a diffusively-coupled network with self-loops.
\end{abstract}

\end{frontmatter}

\section{Introduction}
Network controllability is a topic of active research, with many works in the literature that investigate it from different points of view.
From the viewpoint of the structure of a network, it is common to use only some information about the edges in a network. 
The problem of network controllability for structured networks, which are characterized by the non-zero/zero patterns of edges, was first introduced as \textit{structural controllability} in \cite{lin1974structural}.
A structured network is considered structurally controllable if it can be controlled for almost all choices of non-zero edge weights, 
though it can still be uncontrollable for specific edge weight combinations. This is known as the \textit{generic property} \cite{dion2003generic,park2023fixed} of structural controllability.
To address this limitation, the concept of \textit{strong structural controllability} was introduced in \cite{mayeda1979strong} to guarantee full controllability by considering all choices of edge weights. 
Specifically, a graph is strongly structurally controllable (SSC) if its controllability can be determined based only on the zero/non-zero patterns of edges, independent of the values of the edge weights. 
In the literature, several methods have been proposed to determine the strong structural controllability of structured networks, such as
PMI sequences \cite{yazicioglu2012tight,yaziciouglu2016graph}, zero forcing sets \cite{mousavi2017structural,yaziciouglu2020strong}, 
graph distance \cite{zhang2013upper}, and maximum matching \cite{shabbir2019computation}.
These methods from different perspectives are crucial for a more comprehensive analysis of structured networks, 
which have diverse applications in areas such as targeted networks \cite{monshizadeh2015strong}, social networks  \cite{da2018topology}, and brain networks \cite{gu2015controllability}.

Related to strong structural controllability, the problems of \textit{input addition} and \textit{Minimum Input Problem} (MIP) for control efficiency are important issues. 
These problems have been studied in various approaches, such as the loopy zero forcing set \cite{jia2020unifying}, structural balance \cite{she2018controllability}, 
zero forcing number \cite{monshizadeh2014zero}, and constrained matching \cite{chapman2013strong,trefois2015zero}.
In particular, the authors in \cite{mousavi2017structural,trefois2015zero} proved that the problem of finding minimum inputs for general networks is NP-hard. 
However, it is known that for some special types of small networks, the MIP can be solved in approximate polynomial time using heuristic algorithms,
such as greedy algorithm \cite{olshevsky2014minimal,moothedath2018flow} and independent strongly connected component (iSSC) \cite{guan2021structural}.

In this paper we provide a detailed graph-theoretical interpretation of strong structural controllability, presenting new results as well as reinterpretation of existing conditions. 
Based on these ideas, we offer several methods for designing SSC networks as well as an algorithm with approximate polynomial complexity to solve the MIP problem. 
We conclude this introduction with a guide to the flow of ideas in the paper as well as a summary of our distinct contributions.

\subsection{Research Flow}
This paper introduces several conditions for strong structural controllability based on the new concept of \textit{dedicated} \& \textit{sharing nodes}, 
which is a more detailed concept of \textit{dilation} \cite{lin1974structural} for structural controllability.
For diffusively-coupled networks with self-loops,
we provide a reinterpretation of \textit{Theorem 1} in \cite{mayeda1979strong}, which provides the necessary and sufficient condition for strong structural controllability, from the perspective of the dedicated node.
Based on \cite{ahn2019topological}, we provide conditions for strong structural controllability of basic graph components, i.e.,  paths, cycles, and tree structures.
In particular, we provide an intuitive interpretation by decomposing these basic graph components into paths and employing a \textit{composition rule} for proof. 
The \textit{composition} rule refers to the principle that disjoint controllable components can merge in a manner to ensure that the overall graph meets the controllability conditions.
From this perspective we develop a composition process for pactus type graphs, which we will define in subsequent sections of this paper.
These graphs consist of basic graph components and represent a more generalized concept than cactus, as defined in \cite{menara2018structural}.
Finally, by interpreting the properties of external input nodes from the perspective of dedicated nodes, 
we define the concept of a \textit{component input node}, which is a state node that has the same property as an external input node. 
Based on these two types of input nodes, we then propose an algorithm of approximate polynomial complexity to solve MIP.

\subsection{Contributions}
Note that this paper is an advanced version of \cite{ahn2019topological}.
Although Section~\ref{sec_SSC} in this paper follows a similar logical flow to \cite{ahn2019topological},
there are some significant differences.
First, there is an error in \cite[\textit{Corollary 1}]{ahn2019topological}. 
This paper corrects this error and makes the results flawless.
Furthermore, while \cite{ahn2019topological} is based on the sufficient condition for controllability,
here we provide necessary and sufficient conditions for strong structural controllability.
The contributions of this paper are as follows:
\begin{itemize}
\item
Different from the condition of strong structural controllability introduced in \cite{mayeda1979strong},
we provide a simplified condition for strong structural controllability of undirected graphs of diffusively-coupled networks with self-loops based on the notions of \textit{dedicated} \& \textit{sharing nodes}.

\item 
The paper provides intuitive insights into network controllability by decomposing basic components such as cycles and trees into path graphs and analyzing them based on composition rules.
Furthermore, we expand these findings to more complex graph types, such as pactus, which consist of these basic elements.
This approach highlights the topological control paths within these networks, providing clarity and enhanced understanding.

\item
From the perspective of the \textit{dedicated} \& \textit{sharing nodes},
we present the new notion of a \textit{component input node}, which is a state node that has the same property as an external input node.
Based on this notion,
we devise an algorithm for pactus type graphs that efficiently identifies the minimum external input nodes necessary for strong structural controllability, offering a deeper insight into such structured networks.
\end{itemize}
The paper is organized as follows. In Section~\ref{sec_pre}, the preliminaries and the problems of strong structural controllability are formulated. 
In Section~\ref{sec_SSC}, the conditions for strong structural controllability of basic components and pactus are presented.
Section~\ref{sec_topo_minimum} develops an algorithm for solving MIP for pactus graphs. 
Conclusions are presented in Section~\ref{sec_conc}.

\section{Preliminaries and Problem Formulations} \label{sec_pre}
We will consider undirected networks of diffusively-coupled states $x_i$ and external inputs $u_i$ with self-loops:
\begin{align}
\dot{x}_i = -\sum_{j \in \mathcal{N}_i} a_{ij}(x_i-x_j) + a_{ii}x_i + b_i u_i, \label{eq_1}
\end{align}
where $a_{ij}$ represents the diffusive couplings between states $x_i$ and $x_j$, satisfying $a_{ij}=a_{ji}$, 
$\mathcal{N}_i$ denotes the set of indices of the states $x_j$ for which $a_{ji}\neq 0$, and $b_i$ represents the external input coupling of $x_i$.
This paper assumes that for non-zero values of $a_{ij}$, $a_{ij}$ is positive when $i\neq j$, and is nonzero and negative when $i= j$ as per \cite{yaziciouglu2020strong,shabbir2022computation}.
The Laplacian matrix is defined as $\mathcal{L}=\mathcal{A}-\mathcal{D}$,
where $\mathcal{A}\in\mathbb{R}^{n \times n}$ is the adjacency matrix consisting of diffusive couplings $a_{ij}$ for $i\neq j$,
and $\mathcal{D}=diag(\mathcal{A}\mathds{1}_n)\in\mathbb{R}^{n \times n}$.
Let us now define a self-loop matrix $\mathcal{S}\in\mathbb{R}^{n \times n}$, which is the diagonal matrix consisting of $a_{ii}$ for $i\in\{1,...,n\}$.
Then, the diffusively-coupled network in \eqref{eq_1} can be represented as:
\begin{align}
\dot{x} = \tilde{\mathcal{L}} x + B u, \label{eq_dynamics}
\end{align}
where $x=(x_1,\ldots, x_n)^T\in\mathbb{R}^{n \times 1}$, $u=(u_1, \ldots, u_m)^T\in\mathbb{R}^{m \times 1}$, 
$\tilde{\mathcal{L}}=\mathcal{L}-\mathcal{S} \in \Bbb{R}^{n \times n}$ is the modified Laplacian matrix, and $B=[e_1,e_2,...,e_m]  \in \Bbb{R}^{n \times m} $ is 
the input matrix consisting of linearly independent standard column bases with appropriate dimensions.
Let the diffusively-coupled network matrix corresponding to (\ref{eq_dynamics}) be symbolically written as $T=[\tilde{\mathcal{L}}, B]\in\mathbb{R}^{n \times (n+m)}$.
From the network perspective,
the modified Laplacian matrix $\tilde{\mathcal{L}}$ includes the interactions of $n$ state nodes including each node's own interaction.
The input matrix $B$ includes interactions between $m$ external input nodes and state nodes.
Hence, there are $n+m$ nodes in the network. 

\begin{rem}\label{remark_assum}
The negative weight of the self-loops, $a_{ii}<0$ in \eqref{eq_1}, represents a stabilizing force that prevents uncontrolled growth in the node's state due to self-influence. 
This assumption is commonly used to ensure the stability of the network dynamics \cite{arcak2011diagonal,kaszkurewicz2012matrix}, 
as a positive self-loop term can lead to instability in the system.
\end{rem}


For structured networks, we define a family set of modified Laplacian matrices as $Q(\tilde{\mathcal{L}})$, 
which is determined by the same non-zero/zero patterns as $\mathcal{A}$ and $\mathcal{S}$ in an element-wise fashion. 
Note that since the structured network of \eqref{eq_dynamics} is determined only by the non-zero/zero patterns of $\mathcal{A}$ and $\mathcal{S}$, 
all $\mathcal{L}'\in Q(\tilde{\mathcal{L}})$ have off-diagonal elements with the same non-zero/zero pattern, 
but the diagonal elements of $\mathcal{L}'$ do not necessarily have the same non-zero/zero pattern.
From the viewpoint of control system design, we assume that the input matrix $B$ is fixed since the variations of non-zero elements in $B$ do not affect the controllability.
Thus, the family set of network matrix $T$ is defined as $Q(T) := [Q(\tilde{\mathcal{L}}), B]$.
The network given by \eqref{eq_dynamics} is said to be controllable if its controllability matrix has full row rank \cite{rugh1996linear}.
The controllability matrix corresponding to $T= [\tilde{\mathcal{L}}, B]$ is given by 
$\mathcal{C}_{\tilde{\mathcal{L}}} = [B, \tilde{\mathcal{L}}B, \tilde{\mathcal{L}}^2B, ..., \tilde{\mathcal{L}}^{n-1}B]$.

From a \textit{graph} point of view, 
the network given by \eqref{eq_dynamics} can be represented as a \textit{graph}:
\begin{align}
 \mathcal{G}(T) = (\mathcal{V}, \mathcal{E}),
\end{align}
where $T=[t_{ij}]=[\tilde{\mathcal{L}}, B]$, the node set $\mathcal{V}$ is the union of the set of state nodes and the set of input nodes, 
i.e., $\mathcal{V}=\mathcal{V}^{S}\cup\mathcal{V}^{I}$ satisfying $\mathcal{V}^{S}\cap\mathcal{V}^{I}=\emptyset$, 
and the edge set $\mathcal{E}$ is defined by the interactions between nodes $\mathcal{V}$.
The graph $\mathcal{G}(T)=(\mathcal{V},\mathcal{E})$ consists of a state graph $\mathcal{G}^S$ and 
an interaction graph $\mathcal{G}^I$ as $\mathcal{G}(T) = \mathcal{G}^S \cup  \mathcal{G}^I$, 
where $\mathcal{G}^S$ is the subgraph induced by $\mathcal{V}^S$, and $\mathcal{G}^I$ is the graph representing the interactions between $\mathcal{V}^S$ and $\mathcal{V}^I$. 
That is, $\mathcal{G}^S = (\mathcal{V}^S, \mathcal{E}^S)$ and $\mathcal{G}^I = (\mathcal{V}, \mathcal{E}^{I})$,
the direction of edges in $\mathcal{E}^{S}$ is undirected, while the direction of edges in $\mathcal{E}^{I}$ is directed such that 
$(i,j)\in\mathcal{E}^{I}$ with $i\in\mathcal{V}^{I}$ and $j\in\mathcal{V}^{S}$.
Also, since the input matrix in \eqref{eq_dynamics} consists of the standard column bases, each external input node is connected to only one state node.
Then, the strong structural controllability of a graph corresponding to the diffusively-coupled network given by \eqref{eq_dynamics} is defined as:

\begin{defn}\label{definition_TC}
(Strong structural controllability) A graph $\mathcal{G}(T)$ of \eqref{eq_dynamics} is said to be strongly structurally controllable (SSC)
if all network matrices $T^{'}\in Q(T)$ are controllable.
\end{defn}

We say that the graph $\mathcal{G}(T)$ is \textit{accessible} if there exists a path from $j \in \mathcal{V}^I$ to $i$ for any $i \in \mathcal{V}^S$.
Note that if there is no path from an external input node $i\in\mathcal{V}^{I}$ to a state node $j\in\mathcal{V}^{S}$, then $j$ is not controllable.
Without loss of generality, this paper assumes that a graph $\mathcal{G}(T)$ is accessible, as these are necessary for ensuring network controllability.

\begin{lem}
For a diffusively-coupled network given by \eqref{eq_dynamics}, 
the matrix $\tilde{\mathcal{L}}^{'}$ of a graph $\mathcal{G}(T)$ has full rank for all $\tilde{\mathcal{L}}^{'}\in Q(\tilde{\mathcal{L}})$.
\label{lemma_L_influenceable}
\end{lem}
\begin{pf}
Let us consider the matrix $\tilde{\mathcal{L}}=\mathcal{L}-\mathcal{S}$ in \eqref{eq_dynamics}, 
where $\mathcal{L}\in\mathbb{R}^{n \times n}$ is a Laplacian matrix and $\mathcal{S}\in\mathbb{R}^{n \times n}$ is a diagonal self-loop matrix. 
We assume that the edges between two different state nodes have a positive sign, so every $\mathcal{L}'\in Q(\mathcal{L})$ must be positive semi-definite, 
i.e., $x^T\mathcal{L}'x\geq 0$ for all $x\in\mathbb{R}^n$.
Now, consider the matrix $\mathcal{S}=diag[a_{11},...,a_{nn}]$. Since we assume that every state node in $\mathcal{G}(T)$ has a self-loop with a negative sign,
 i.e., $a_{ii}<0$ for all $i\in{1,...,n}$, every $\mathcal{S}$ must be negative definite, implying that $x^T(-\mathcal{S})x>0$ for all $x\in\mathbb{R}^n\setminus\{0\}$.
Therefore, we obtain $x^T\tilde{\mathcal{L}}'x=x^T(\mathcal{L}'-\mathcal{S})x>0$ for all $x\in\mathbb{R}^n\setminus\{0\}$, 
which means that $\tilde{\mathcal{L}}'$ is positive definite. Hence, every $\tilde{\mathcal{L}}' \in Q(\tilde{\mathcal{L}})$ has full rank.
\hfill $\square$
\end{pf}

From the above lemma, it follows that every network matrix $T'\in Q(T)$ for the diffusively-coupled network with self-loops given by \eqref{eq_dynamics} also has a full row rank.

Next, let the neighboring set of the set $\alpha \subseteq \mathcal{V}^S$ be denoted as $\mathcal{N}(\alpha)=\bigcup_{j\in\alpha} \mathcal{N}_j$, which is the union set of neighboring set of $j$ satisfying $j\in\alpha$.
From \cite[\textit{Lemma 3}]{mayeda1979strong}, the conditions provided in \cite[\textit{Theorem 1}]{mayeda1979strong} provided can be simplified for diffusively-coupled networks with self-loops as:


\begin{thm} \cite{mayeda1979strong}
The graph $\mathcal{G}(T)$ is SSC if and only if for all $\alpha \subseteq \mathcal{V}^S$, 
there exists $i\in\mathcal{N}(\alpha)$$\setminus$$\alpha$ having exactly one edge $(i,j) \in \mathcal{E}$ with $j \in \alpha$ in $\mathcal{G}(T)$.\footnote{Note that in this theorem and in our results below, it is necessary to consider all $\alpha \subseteq \mathcal{V}^S$. The number of possible sets of $\alpha$ is $2^{\vert V^S \vert}-1$.}
\label{theorem_Tsatsomeros}
\end{thm}

For a more detailed analysis, we classify the node $i\in\mathcal{N}(\alpha)$$\setminus$$\alpha$ as \textit{dedicated nodes} and \textit{sharing nodes} for $\alpha\subseteq\mathcal{V}^{S}$.

\begin{defn}\label{definition_candidate_node}
(Dedicated and Sharing nodes) 
A node $i\in\mathcal{N}(\alpha)$$\setminus$$\alpha$ is a dedicated node of $\alpha$ if it satisfies $|\mathcal{N}_{i}\cap\alpha|=1$,
or is a sharing node of $\alpha$ if it satisfies $|\mathcal{N}_{i}\cap\alpha|>1$.
\end{defn}

For an arbitrary $\alpha$ satisfying $\alpha \subseteq \mathcal{V}^S$,
if a node $i\in\mathcal{N}(\alpha)$$\setminus$$\alpha$ has exactly one edge connected to $\alpha$,
then the node $i$ is called a dedicated node of $\alpha$.
This statement is equivalent to the cardinality condition of $|\mathcal{N}_{i}\cap\alpha|=1$.
On the other hand, a node $i\in\mathcal{N}(\alpha)$$\setminus$$\alpha$ is called a sharing node if the node $i$ has more than one edge connected to $\alpha$.
Using the concepts of dedicated and sharing nodes, a reinterpretation of \textit{Theorem~\ref{theorem_Tsatsomeros}} is that a 
a graph $\mathcal{G}(T)$ is SSC if and only if the set $\mathcal{N}(\alpha)$$\setminus$$\alpha$ has at least one dedicated node
for all $\alpha\subseteq\mathcal{V}^{S}$. This is stated formally below in \textit{Corollary~\ref{corollary_theorem1}}.

For example, consider the graph $\mathcal{G}(T)$ depicted in Fig.~\ref{network_ex_tree}(a). 
It is shown that $\mathcal{V}^{S}=\{1,2,3,4,5\}$ and $\mathcal{V}^{I}=\{u_1,u_2\}$.
Let $\alpha \subseteq \mathcal{V}^S$ be $\alpha=\{1,3\}$.
Then, we obtain $\mathcal{N}(\alpha)$$\setminus$$\alpha=\{2,u_1\}$.
Now, we need to check whether the set $\mathcal{N}(\alpha)$$\setminus$$\alpha$ has a dedicated node or not.
For the node $2\in$ $\mathcal{N}(\alpha)$$\setminus$$\alpha$, we obtain $\mathcal{N}_{2}=\{1,3,4\}$.
Then, it follows that $|\mathcal{N}_{2}\cap\alpha|=2$, which means that the node $2$ is a sharing node
since those have two edges $(2,1),(2,3)\in\mathcal{E}$ connected to $1,3\in\alpha$, respectively. 
For the external input node $u_1\in$ $\mathcal{N}(\alpha)$$\setminus$$\alpha$, we obtain $\mathcal{N}_{u_1}=\{3\}$.
Then, we obtain $|\mathcal{N}_{u_1}\cap\alpha|=1$, which means that $u_1$ is a dedicated node of $\alpha$
because the node $u_1$ has exactly one edge $(u_1,3)\in\mathcal{E}$ connected to $3\in\alpha$.
Hence, in case of $\alpha=\{1,3\}$, there exists a dedicated node $u_1$.
In the same way, if the set $\mathcal{N}(\alpha)$$\setminus$$\alpha$ has at least one dedicated node 
for all $\alpha\subseteq\mathcal{V}^{S}$, the graph $\mathcal{G}(T)$ is determined to be a SSC graph.

\begin{rem}
Note that the concept of \textit{dedicated} \& \textit{sharing node} is an extended concept of \textit{dilation} \cite{lin1974structural} for structural controllability.
More precisely, the authors in \cite{lin1974structural} show that the condition of structural controllability is the absence of \textit{dilation} in a graph,
but this condition is not sufficient for strong structural controllability. \textit{Dedicated} \& \textit{sharing nodes} are a more detailed concept for determining the strong structural controllability.
\end{rem}

\section{Strongly Structurally Controllable Graphs} \label{sec_SSC}
In this section, we provide several necessary and sufficient conditions for strong structural controllability of basic graph components such as paths, cycles, and tree structures. 
While there is a wealth of research on determining controllability conditions for these graph types using concepts, 
such as the \textit{zero forcing set} \cite{mousavi2019strong} and \textit{structural balance} \cite{she2018controllability}, 
this paper offers a reinterpretation of these existing conditions.
This is driven by the belief that alternative approaches to the same results can provide deeper insights and different perspectives on the problem at hand.
We specifically analyze the controllability of basic components based on the notions of dedicated and sharing nodes, employing composition rules to deepen the analysis.
With the notion of a \textit{dedicated node},
\textit{Theorem~\ref{theorem_Tsatsomeros}} can be simplified as:
\begin{cor}\label{corollary_theorem1} 
The graph $\mathcal{G}(T)$ is SSC
if and only if there exists at least one dedicated node in $\mathcal{N}(\alpha)$$\setminus$$\alpha$
for all $\alpha \subseteq \mathcal{V}^{S}$.
\end{cor}

The above \textit{Corollary~\ref{corollary_theorem1}} provides the necessary and sufficient condition for strong structural controllability from the perspective of dedicated nodes.
For further analysis, the following proposition provides the condition for strong structural controllability of the path graph, a basic component, as presented in \cite{mayeda1979strong}.
\begin{prop}\label{proposition_path}
\cite{mayeda1979strong}
Let us consider a path state graph $\mathcal{G}^{S}=(\mathcal{V}^{S},\mathcal{E}^{S})$
with $\mathcal{G}^{I}=(\mathcal{V}^{I},\mathcal{E}^{I})$. 
The graph $\mathcal{G}(T)=\mathcal{G}^{S}\cup\mathcal{G}^{I}$ is SSC
if and only if there exists an external input node connected to a terminal state node\footnote{In a graph $\mathcal{G}(T)=(\mathcal{V},\mathcal{E})$, we say that a state node $k\in\mathcal{V}^{S}$ is a terminal state node if it satisfies $|\mathcal{N}_k|=1$ in $\mathcal{G}^{S}$.}.
\end{prop}
\begin{figure}[]
\centering
\subfigure[]{
\begin{tikzpicture}[scale=0.8]
\node[place, black] (node1) at (-2,0) [label=below:\scriptsize$1$] {};
\node[place, black] (node2) at (-1,1) [label=below:\scriptsize$2$] {};
\node[place, black] (node3) at (0,2) [label=below:\scriptsize$3$] {};
\node[place, black] (node4) at (0,0.5) [label=below:\scriptsize$4$] {};
\node[place, black] (node5) at (1,0) [label=below:\scriptsize$5$] {};

\node[place, circle] (node6) at (0.7,2.7) [label=above:\scriptsize$u_1$] {};
\node[place, circle] (node7) at (1.7,0.7) [label=above:\scriptsize$u_2$] {};

\draw (node1) [line width=0.5pt] -- node [left] {} (node2);
\draw (node2) [line width=0.5pt] -- node [right] {} (node3);
\draw (node2) [line width=0.5pt] -- node [left] {} (node5);
\draw (node6) [-latex, line width=0.5pt] -- node [right] {} (node3);
\draw (node7) [-latex, line width=0.5pt] -- node [right] {} (node5);
\end{tikzpicture}
}
\subfigure[]{
\begin{tikzpicture}[scale=0.8]
\node[place, black] (node1) at (-2,0) [label=below:\scriptsize$1$] {};
\node[place, black] (node2) at (-1,1) [label=below:\scriptsize$2$] {};
\node[place, black] (node3) at (0,2) [label=below:\scriptsize$3$] {};
\node[place, black] (node4) at (0,0.5) [label=below:\scriptsize$4$] {};
\node[place, black] (node5) at (1,0) [label=below:\scriptsize$5$] {};

\node[place, circle] (node6) at (0.7,1.2) [label=above:\scriptsize$u_1$] {};
\node[place, circle] (node7) at (1.7,0.7) [label=above:\scriptsize$u_2$] {};

\draw (node1) [line width=0.5pt] -- node [left] {} (node2);
\draw (node2) [line width=0.5pt] -- node [right] {} (node3);
\draw (node2) [line width=0.5pt] -- node [left] {} (node5);
\draw (node6) [-latex, line width=0.5pt] -- node [right] {} (node4);
\draw (node7) [-latex, line width=0.5pt] -- node [right] {} (node5);
\end{tikzpicture}
}

\caption{Tree state graphs with two external input nodes.}
\label{network_ex_tree}
\centering
\subfigure[]{
\begin{tikzpicture}[scale=0.8]
\node[place, black] (node1) at (-2,0) [label=above:\scriptsize$1$] {};
\node[place, black] (node2) at (-1,1) [label=above:\scriptsize$2$] {};
\node[place, black] (node3) at (0,-1) [label=below:\scriptsize$3$] {};
\node[place, black] (node4) at (-1,-1) [label=below:\scriptsize$4$] {};

\node[place, circle] (node5) at (0,1) [label=below:\scriptsize$u_1$] {};
\node[place, circle] (node6) at (1,-1) [label=below:\scriptsize$u_2$] {};

\draw (node1) [line width=0.5pt] -- node [left] {} (node2);
\draw (node2) [line width=0.5pt] -- node [right] {} (node3);
\draw (node3) [line width=0.5pt] -- node [right] {} (node4);
\draw (node1) [line width=0.5pt] -- node [left] {} (node4);
\draw (node5) [-latex, line width=0.5pt] -- node [right] {} (node2);
\draw (node6) [-latex, line width=0.5pt] -- node [right] {} (node3);
\end{tikzpicture}
}
\subfigure[]{
\begin{tikzpicture}[scale=0.8]
\node[place, black] (node1) at (-2,0) [label=below:\scriptsize$1$] {};
\node[place, black] (node2) at (-1,1) [label=above:\scriptsize$2$] {};
\node[place, black] (node3) at (0,-1) [label=below:\scriptsize$3$] {};
\node[place, black] (node4) at (-1,-1) [label=below:\scriptsize$4$] {};

\node[place, circle] (node5) at (-3,0) [label=below:\scriptsize$u_1$] {};
\node[place, circle] (node6) at (1,0) [label=left:\scriptsize$u_2$] {};

\draw (node1) [line width=0.5pt] -- node [left] {} (node2);
\draw (node2) [line width=0.5pt] -- node [right] {} (node3);
\draw (node3) [line width=0.5pt] -- node [right] {} (node4);
\draw (node1) [line width=0.5pt] -- node [left] {} (node4);
\draw (node5) [-latex, line width=0.5pt] -- node [right] {} (node1);
\draw (node6) [-latex, line width=0.5pt] -- node [right] {} (node3);

\end{tikzpicture}
}
\caption{Cycle state graphs with two external input nodes.}
\label{network_ex_cycle}

\end{figure}
For further analysis of a larger graph, 
we define a bridge graph $\mathcal{G}_{ij}^{S}$, 
which connects two disjoint graphs $\mathcal{G}_i$ and $\mathcal{G}_j$.

\begin{defn}\label{definition_Bridge_graph} 
(Bridge graph) 
A bridge graph is a state graph defined as $\mathcal{G}_{ij}^{S}=(\mathcal{V}_{ij}^{S},\mathcal{E}_{ij}^{S}$),
which connects two disjoint state graphs $\mathcal{G}^{S}_{i}$ and $\mathcal{G}^{S}_{j}$.
If nodes $k\in\mathcal{V}^{S}_{i}$ and $l\in\mathcal{V}^{S}_{j}$ are connected by an edge $(k,l)$, then $(k,l)\in\mathcal{E}_{ij}^{S}$ and $k,l\in\mathcal{V}_{ij}^{S}$.
We assume that a state node $k\in\mathcal{V}^{S}_{i}$ is connected to a state node $l\in\mathcal{V}^{S}_{j}$ by a one-to-one.
Hence, $|\mathcal{E}_{ij}^{S}|$ satisfies the following boundary condition.
\begin{align}\label{bridge_edge_boundary}
1\le |\mathcal{E}_{ij}^{S}|\le min(|\mathcal{V}_{i}^{S}|,|\mathcal{V}_{j}^{S}|).
\end{align}
\end{defn} 

For convenience, we first define the concept of neighbor components for a given component $\mathcal{G}_i^S$ in a graph. 
The set of neighbor components, denoted as $\mathcal{N}_{\mathcal{G}_i^S}$, 
includes the indices $j$ for which there exists a non-empty bridge graph $\mathcal{G}_{ij}^S$ connecting $\mathcal{G}_i^S$ to $\mathcal{G}_j^S$. 
This is formally expressed as $\mathcal{N}_{\mathcal{G}_i^S} = \{ j \mid \mathcal{G}_{ij}^S \neq \emptyset \}$, 
where each index $j$ in this set represents a component directly connected to $\mathcal{G}_i^S$ by at least bridge edge. 
With the above notation, we say that the graph $\mathcal{G}(T)=(\mathcal{V},\mathcal{E})$ is induced by 
$m$-disjoint components $\mathcal{G}_{i}=(\mathcal{V}_{i},\mathcal{E}_{i})$ with bridge graph edges $\mathcal{E}_{ij}^{S}$
if $\mathcal{V}={\bigcup}^{m}_{i=1}\mathcal{V}_{i}$ and $\mathcal{E}={\bigcup}^{m}_{i=1}{\bigcup}_{j\in\mathcal{N}_{\mathcal{G}_i^S}}(\mathcal{E}_{i}    \cup    \mathcal{E}_{ij}^{S})$.
Based on the concept above of induction, we now proceed to decompose larger graphs into the basic components, i.e., path graphs, 
and explore their strong structural controllability conditions.

\begin{thm}\label{theorem_tree}
Consider a tree state graph $\mathcal{G}^{S}=(\mathcal{V}^{S},\mathcal{E}^{S})$
with $\mathcal{G}^{I}=(\mathcal{V}^{I},\mathcal{E}^{I})$ satisfying $|\mathcal{V}^{I}|=m\ge2$. 
Then, the graph $\mathcal{G}(T)=\mathcal{G}^{S}\cup\mathcal{G}^{I}$ can be induced by $m$-disjoint components $\mathcal{G}_{i}$ with $\mathcal{E}_{ij}^{S}$ 
satisfying $|\mathcal{E}_{ij}^{S}|=1$ and $|\mathcal{V}^{I}_{i}|=1$ for $i\in \{1,...,m\}$ and $j\in\mathcal{N}_{\mathcal{G}_i^S}$.
The graph $\mathcal{G}(T)=\mathcal{G}^{S}\cup\mathcal{G}^{I}$ is SSC
if and only if each disjoint component $\mathcal{G}_{i}, i\in\{1,...,m\}$ is SSC path graph satisfying \textit{Proposition~\ref{proposition_path}}.
\end{thm}
\begin{pf}
For \textit{if} condition,
let us assume that each disjoint component $\mathcal{G}_{i},i\in\{1,...,m\}$ satisfies \textit{Proposition~\ref{proposition_path}}.
Then, since each $\mathcal{G}_{i}$ is SSC, 
the set $\mathcal{N}(\alpha_{i})$$\setminus$$\alpha_{i}$ contains at least one dedicated node
for all $\alpha_{i}\subseteq\mathcal{V}_{i}^{S}$.
Now, consider the merged graph $\mathcal{G}(T)$ with the bridge edge $(k,l)\in\mathcal{E}_{ij}^{S}$, 
where $k\in\mathcal{V}_{i}^{S}$ and $l\in\mathcal{V}_{j}^{S}$ for $i\in \{1,...,m\}$ and $j\in\mathcal{N}_{\mathcal{G}_i^S}$.
For the merged graph $\mathcal{G}(T)$ to be SSC, we only need to consider the existence of dedicated nodes in $\mathcal{N}(\alpha)$$\setminus$$\alpha$
when $\alpha\subseteq\mathcal{V}^{S}$ contains at least one bridge node. 
Because after adding the bridge edge, each set of in-neighboring nodes of a node in $\mathcal{V}^{S}$ 
remains unchanged except for $\mathcal{N}_{k}$ and $\mathcal{N}_{l}$ in $\mathcal{G}(T)$.
However, since each disjoint component $\mathcal{G}_{i},i\in\{1,...,m\}$ is SSC, even if the bridge node $k\in\mathcal{V}_{i}^{S}$ or $l\in\mathcal{V}_{j}^{S}$ belongs to $\alpha\subseteq\mathcal{V}^{S}$, the set $\mathcal{N}(\alpha)$$\setminus$$\alpha$ still has at least one dedicated node in $\mathcal{V}_{i}$ or in $\mathcal{V}_{j}$.
It follows that the existence of at least one dedicated node in $\mathcal{N}(\alpha)$$\setminus$$\alpha$ is independent of 
the bridge edge $(k,l)\in\mathcal{E}_{ij}^{S}$ satisfying $|\mathcal{E}_{ij}^{S}|=1$ for all $\alpha\subseteq\mathcal{V}^{S}$.

For \textit{only if} condition, in the merged graph $\mathcal{G}(T)$,
let us suppose that a disjoint component $\mathcal{G}_{q}=\mathcal{G}_{q}^{S}\cup\mathcal{G}_{q}^{I},q\in\{1,...,m\}$ does not satisfy \textit{Proposition~\ref{proposition_path}}.
Since the graph $\mathcal{G}(T)$ is a tree graph, there always exists a state node $i\in\mathcal{V}^{S}$, which has out-degree 3.
Then, there always exists a case without a dedicated node in $\mathcal{N}(\alpha)$$\setminus$$\alpha$
when $\alpha=\mathcal{V}^{S}$$\setminus$$\{i\}$.
\hfill $\square$
\end{pf}
As an example of \textit{Theorem~\ref{theorem_tree}}, the graph $\mathcal{G}(T)$ depicted in Fig.~\ref{network_ex_tree}(a) shows a tree state graph $\mathcal{G}^{S}$ with two external input nodes $u_1,u_2\in\mathcal{V}^{I}$. In this case, the graph $\mathcal{G}(T)$ can be induced by $2$-disjoint path graphs $\mathcal{G}_{1}$ and $\mathcal{G}_{2}$ with a bridge edge $(2,4)\in\mathcal{E}_{12}$, i.e., $\mathcal{G}_{1}:u_1\rightarrow 3 \leftrightarrow 2 \leftrightarrow 1$ and $\mathcal{G}_{2}:u_2\rightarrow 5 \leftrightarrow 4$, where the symbols $\rightarrow$ and $\leftrightarrow$ are used to denote directions of the connection between nodes. It follows from \textit{Theorem~\ref{theorem_tree}} that the merged graph $\mathcal{G}(T)$ is SSC since each disjoint component $\mathcal{G}_{1}$ and $\mathcal{G}_{2}$ satisfies \textit{Proposition~\ref{proposition_path}}. However, the graph $\mathcal{G}(T)$ depicted in Fig.~\ref{network_ex_tree}(b) can not be induced by 2-disjoint path graphs. In this case, when $\alpha=\{ 1,3\}$, we obtain $\mathcal{N}(\alpha)$$\setminus$$\alpha=\{2\}$. But the node $2$ is a sharing node satisfying $|\mathcal{N}_{2}\cap\alpha|>1$, which is connected to the nodes $1,3\in\alpha$, thus, the graph $\mathcal{G}(T)$ in Fig.~\ref{network_ex_tree}(b) is not SSC.
With the result of \textit{Theorem~\ref{theorem_tree}}, the following \textit{Corollary~\ref{corollary_1}} can be obtained:

\begin{cor}\label{corollary_1}
Let two components $\mathcal{G}_{i}$ and $\mathcal{G}_{j}$ be SSC, respectively.
If there exists a bridge graph $\mathcal{G}_{ij}^{S}$ satisfying $|\mathcal{E}_{ij}^{S}|=1$,
then the merged graph $\mathcal{G}(T)=\mathcal{G}_{i}\cup\mathcal{G}_{ij}^{S}\cup\mathcal{G}_{j}$ is SSC, regardless of the location of the bridge edge.
\end{cor}

The above corollary means that the existence of one bridge edge connecting two disjoint SSC graphs
is independent of the strong structural controllability of the merged graph.

\begin{thm}\label{theorem_cycle}
Consider a cycle state graph $\mathcal{G}^{S}=(\mathcal{V}^{S},\mathcal{E}^{S})$
with $\mathcal{G}^{I}=(\mathcal{V}^{I},\mathcal{E}^{I})$ satisfying $|\mathcal{V}^{I}|=2$. 
The graph $\mathcal{G}(T)=\mathcal{G}^{S}\cup\mathcal{G}^{I}$ is SSC
if and only if there exists an edge $(k,l)\in\mathcal{E}$ with $ k,l \in\mathcal{N}(\mathcal{V}^{I})$.
\end{thm}

\begin{pf} 
Let us consider that a graph $\mathcal{G}(T)=\mathcal{G}^{S}\cup\mathcal{G}^{I}$ consists of a cycle state graph $\mathcal{G}^{S}=(\mathcal{V}^{S},\mathcal{E}^{S})$ 
with $\mathcal{G}^{I}=(\mathcal{V}^{I},\mathcal{E}^{I})$ satisfying $|\mathcal{V}^{I}|=2$.
Then, the graph $\mathcal{G}(T)$ can be induced by $2$-disjoint components 
$\mathcal{G}_{1}$ and $\mathcal{G}_{2}$ with $\mathcal{E}_{12}^{S}$ satisfying $|\mathcal{E}_{12}^{S}|=2$.
Also, each disjoint component $\mathcal{G}_{1}$ and $\mathcal{G}_{2}$ satisfies \textit{Proposition~\ref{proposition_path}}.
Thus, the sets $\mathcal{N}(\alpha_{1})$$\setminus$$\alpha_{1}$ and $\mathcal{N}(\alpha_{2})$$\setminus$$\alpha_{2}$ have at least one dedicated node
for all $\alpha_{1}\subseteq\mathcal{V}_{1}^{S}$ and $\alpha_{2}\subseteq\mathcal{V}_{2}^{S}$, respectively.
For the merged graph $\mathcal{G}(T)$ to be SSC,
we only need to consider the existence of dedicated nodes in $\mathcal{N}(\alpha)$$\setminus$$\alpha$
when $\alpha\subseteq\mathcal{V}^{S}$ contains at least one bridge node.
Because each set of in-neighboring nodes of a node in $\mathcal{V}^{S}$ remains unchanged except for the nodes in $\mathcal{V}_{12}^{S}$.
Now, start from $\mathcal{G}_{1}\cup\mathcal{G}_{2}$, we gradually add two bridge edges step-by-step for check the condition of \textit{Corollary~\ref{corollary_theorem1}}.
Let the bridge edges be $\{(k_1,l_1),(k_2,l_2)\}\in\mathcal{E}^{S}_{12}$, 
where $k_1,k_2\in\mathcal{V}_{1}^{S}$ and $l_1,l_2\in\mathcal{V}_{2}^{S}$ satisfying $k_1,l_1\notin\mathcal{N}(\mathcal{V}^{I})$.

For \textit{if} condition,
consider a merged graph $\mathcal{G}_{1}\cup\mathcal{G}_{2}$ with the bridge edge $(k_1,l_1)\in\mathcal{E}^{S}_{12}$.
It follows from \textit{Corollary~\ref{corollary_1}} that if each $\mathcal{G}_{1}$ and $\mathcal{G}_{2}$ is SSC, 
the merged graph $\mathcal{G}_{1}\cup\mathcal{G}_{2}$ with a bridge edge $(k_1,l_1)\in\mathcal{E}_{12}^{S}$ is SSC.
For the other bridge edge $(k_2,l_2)\in\mathcal{E}^{S}_{12}$,
suppose that the bridge nodes $k_2,l_2$ satisfy $k_2,l_2\in\mathcal{N}(\mathcal{V}^{I})$.
Then, each set of in-neighboring nodes of nodes $k_2$ and $l_2$ always includes each other, i.e., $k_2\in\mathcal{N}_{l_2}$ and $l_2\in\mathcal{N}_{k_2}$.
Hence, if $\alpha\subseteq\mathcal{V}^{S}$ contains $k_2$ or $l_2$,
the set $\mathcal{N}(\alpha)$$\setminus$$\alpha$ always contains at least one of the nodes $k_2$ and $l_2$ with $k_2,l_2\in\mathcal{N}(\mathcal{V}^{I})$.
However, if $\alpha$ contains a node in $\mathcal{N}(\mathcal{V}^{I})$, there always exists at least one dedicated node in $\mathcal{N}(\alpha)$$\setminus$$\alpha$,
and this dedicated node is always an external input node in $\mathcal{V}^{I}$.
Therefore, the graph $\mathcal{G}(T)$ is SSC.
For \textit{only if} condition,
let us suppose that $(k_2,l_2)\notin\mathcal{E}$ with $k_2,l_2\in\mathcal{N}(\mathcal{V}^{I})$.
In this case, when choosing $\alpha=\mathcal{V}^{S}$,
we obtain $\mathcal{N}(\alpha)$$\setminus$$\alpha=\{k_2,l_2\}$.
However, both nodes $k_2$ and $l_2$ are sharing nodes satisfying $|\mathcal{N}_{k_2}\cap\alpha|=|\mathcal{N}_{l_2}\cap\alpha|=2$. 
Hence, from \textit{Corollary~\ref{corollary_theorem1}}, the graph $\mathcal{G}(T)$ is not SSC.
\hfill $\square$
\end{pf}
%

For example, Fig.~\ref{network_ex_cycle}(a) shows a cycle state graph $\mathcal{G}^{S}$ with $\mathcal{N}(\mathcal{V}^{I})=\{ 2,3\}$ and there exists an edge $(2,3)\in\mathcal{E}$. 
According to \textit{Theorem~\ref{theorem_cycle}}, the graph $\mathcal{G}(T)=\mathcal{G}^{S}\cup\mathcal{G}^{I}$ is SSC. 
However, the graph $\mathcal{G}(T)$ depicted in Fig.~\ref{network_ex_cycle}(b) shows $\mathcal{N}(\mathcal{V}^{I})=\{ 1,3\}$ 
and there is no edge between the nodes $1,3\in\mathcal{N}(\mathcal{V}^{I})$, i.e., $(1,3)\notin\mathcal{E}$. 
In this case, when choosing $\alpha=\mathcal{V}^{S}$$\setminus$$\mathcal{N}(\mathcal{V}^{I})=\{ 2,4\}$, we obtain $\mathcal{N}(\alpha)$$\setminus$$\alpha=\{1,3\}$.
It is clear that nodes $1,3\in\mathcal{N}(\alpha)$$\setminus$$\alpha$, are sharing nodes satisfying $|\mathcal{N}_{1}\cap\alpha|=|\mathcal{N}_{3}\cap\alpha|=2$. Hence, according to \textit{Corollary~\ref{corollary_theorem1}}, the graph $\mathcal{G}(T)$ in Fig.~\ref{network_ex_cycle}(b) is not SSC. Note that if $\mathcal{G}^{S}$ is a cycle, two properly located external input nodes are sufficient for the graph $\mathcal{G}(T)$ to be SSC, i.e., the minimum number of external input nodes for the strong structural controllability of $\mathcal{G}(T)$ is $2$. 

With the above proof, \textit{Theorem~\ref{theorem_cycle}} can be generalized as:

\begin{thm}\label{theorem_2}
Let two disjoint components $\mathcal{G}_{i}$ and $\mathcal{G}_{j}$ be SSC with path state graphs $\mathcal{G}^{S}_{i}$ and $\mathcal{G}^{S}_{j}$, respectively.
If there exists a bridge edge $(k,l)\in\mathcal{E}_{ij}^{S}$ satisfying $k,l\in\mathcal{N}(\mathcal{V}^{I})$,
the merged graph $\mathcal{G}_{i}\cup\mathcal{G}_{ij}^{S}\cup\mathcal{G}_{j}$ is SSC, regardless of the existence of an additional bridge edge in $\mathcal{E}_{ij}^{S}$.
\end{thm}

The above \textit{Theorem~\ref{theorem_2}} contains the condition of strong structural controllability 
for a graph $\mathcal{G}(T)=\mathcal{G}^{S}\cup\mathcal{G}^{I}$ when $\mathcal{G}^{S}$ is a cycle.
Thus, \textit{Theorem~\ref{theorem_cycle}} is a special case of \textit{Theorem~\ref{theorem_2}}.
To extend our results to larger graphs, we define a pactus as a structure consisting of disjoint components interconnected by bridge graphs.

\begin{defn}\label{definition_sym_pactus} 
(Pactus) A pactus is a connected graph 
defined as $\mathcal{G}^{S}={\bigcup}^{m}_{i=1}{\bigcup}_{j\in\mathcal{N}_{\mathcal{G}_i^S}}(\mathcal{G}_{i}^{S}    \cup    \mathcal{G}_{ij}^{S})$.
A pactus satisfies the following properties.
\begin{enumerate} 
\item $\mathcal{G}^{S}$ is induced by $m$-disjoint components $\mathcal{G}_{i}^{S}$ with $\mathcal{E}_{ij}^{S}$.
\item Each $\mathcal{G}_{i}^{S}, i\in\{1,...,m\}$, is either a path or a cycle.
\item If $j\in\mathcal{N}_{\mathcal{G}_i^S}$, $\mathcal{G}_{i}^{S}$ and $\mathcal{G}_{j}^{S}$ are connected by at least one bridge edge $(k,l)\in\mathcal{E}_{ij}^{S}$,
where $k\in\mathcal{V}^{S}_{i}$, $l\in\mathcal{V}^{S}_{j}$.
\end{enumerate}
\end{defn}
Note that the pactus is a more generalized concept than the cactus defined in \cite{menara2018structural}.
It means that the cactus is a special case of the pactus.
For example, the bridge edges between two disjoint components in pactus may be several under the boundary condition given by \eqref{bridge_edge_boundary}, while the cactus has only one.
Based on the aforementioned lemmas, the following theorem can be established.

\begin{figure}[]
\centering
\begin{tikzpicture}[scale=0.5]

\node[] at (3.5,4.7) {\scriptsize$\mathcal{G}_1$};
\node[] at (7,6) {\scriptsize$\mathcal{G}_2$};
\node[] at (10.4,3) {\scriptsize$\mathcal{G}_3$};
\node[] at (14,4.5) {\scriptsize$\mathcal{G}_4$};

\node[red] at (6.1,3) {\scriptsize$\mathcal{G}_{12}^{S}$};
\node[red] at (10,5) {\scriptsize$\mathcal{G}_{23}^{S}$};
\node[red] at (11.9,3) {\scriptsize$\mathcal{G}_{34}^{S}$};

\node[place, black] (node1) at (1,5) [label=below:\scriptsize$1$] {};
\node[place, black] (node2) at (3,4) [label=below:\scriptsize$2$] {};
\node[place, red] (node3) at (5,3) [label=below:\scriptsize$3$] {};
\node[place, red] (node4) at (7,4) [label=below:\scriptsize$4$] {};
\node[place, black] (node5) at (5,5) [label=below:\scriptsize$5$] {};
\node[place, black] (node6) at (5,7) [label=above:\scriptsize$6$] {};
\node[place, black] (node7) at (7,8) [label=above:\scriptsize$7$] {};
\node[place, black] (node8) at (9,7) [label=above:\scriptsize$8$] {};
\node[place, red] (node9) at (9,5) [label=below:\scriptsize$9$] {};
\node[place, black] (node10) at (9,3) [label=below:\scriptsize$10$] {};
\node[place, red] (node11) at (11,2) [label=below:\scriptsize$11$] {};
\node[place, red] (node12) at (11,4) [label=above:\scriptsize$12$] {};
\node[place, black] (node13) at (13,5) [label=above:\scriptsize$13$] {};
\node[place, black] (node14) at (15,6) [label=above:\scriptsize$14$] {};
\node[place, black] (node15) at (15,4) [label=below:\scriptsize$15$] {};
\node[place, red] (node16) at (13,3) [label=below:\scriptsize$16$] {};

\draw (node1) [line width=0.5pt] -- node [left] {} (node2);
\draw (node2) [line width=0.5pt] -- node [left] {} (node3);
\draw (node3) [red,dashed,line width=0.5pt] -- node [below] {} (node4);
\draw (node4) [line width=0.5pt] -- node [left] {} (node5);
\draw (node5) [line width=0.5pt] -- node [left] {} (node6);
\draw (node6) [line width=0.5pt] -- node [left] {} (node7);
\draw (node7) [line width=0.5pt] -- node [left] {} (node8);
\draw (node8) [line width=0.5pt] -- node [left] {} (node9);
\draw (node4) [line width=0.5pt] -- node [left] {} (node9);
\draw (node9) [red,dashed,line width=0.5pt] -- node [below] {} (node12);
\draw (node10) [line width=0.5pt] -- node [left] {} (node11);
\draw (node10) [line width=0.5pt] -- node [left] {} (node12);
\draw (node11) [line width=0.5pt] -- node [left] {} (node12);
\draw (node11) [red,dashed,line width=0.5pt] -- node [below] {} (node16);
\draw (node13) [line width=0.5pt] -- node [left] {} (node14);
\draw (node14) [line width=0.5pt] -- node [left] {} (node15);
\draw (node15) [line width=0.5pt] -- node [left] {} (node16);
\draw (node13) [line width=0.5pt] -- node [left] {} (node16);

\node[place, circle] (node17) at (-0.5,4) [label=below:\scriptsize$u_{1}$] {}; 
\node[place, circle] (node19) at (3.5,6) [label=above:\scriptsize$u_{2}$] {}; 
\node[place, circle] (node20) at (3.5,8) [label=above:\scriptsize$u_{3}$] {}; 
\node[place, circle] (node21) at (7.5,2) [label=below:\scriptsize$u_{4}$] {}; 
\node[place, circle] (node22) at (9.5,1) [label=below:\scriptsize$u_{5}$] {}; 
\node[place, circle] (node23) at (11.5,6) [label=above:\scriptsize$u_{6}$] {}; 
\node[place, circle] (node24) at (13.5,7) [label=above:\scriptsize$u_{7}$] {}; 
\draw (node17) [-latex, line width=0.5pt] -- node [left] {} (node1);
\draw (node19) [-latex, line width=0.5pt] -- node [right] {} (node5);
\draw (node20) [-latex, line width=0.5pt] -- node [right] {} (node6);
\draw (node21) [-latex, line width=0.5pt] -- node [right] {} (node10);
\draw (node22) [-latex, line width=0.5pt] -- node [right] {} (node11);
\draw (node23) [-latex, line width=0.5pt] -- node [right] {} (node13);
\draw (node24) [-latex, line width=0.5pt] -- node [right] {} (node14);

\end{tikzpicture}
\caption{A SSC graph $\mathcal{G}(T)=\mathcal{G}^{S}\cup\mathcal{G}^{I}$ with a pactus 
$\mathcal{G}^{S}={\bigcup}^{4}_{i=1}{\bigcup}_{j\in\mathcal{N}_{\mathcal{G}_i^S}}(\mathcal{G}_{i}^{S}    \cup    \mathcal{G}_{ij}^{S})$ satisfying $|\mathcal{E}_{ij}^{S}|=1$.}
\label{network_ex_th2}
\end{figure}
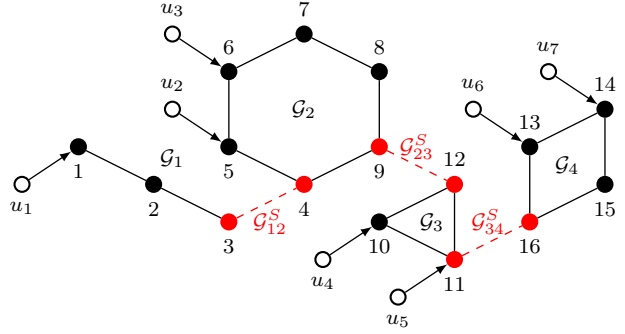

\begin{thm} \label{theorem_sym_pactus_sc_condition}
Let us consider a pactus $\mathcal{G}^{S}$ and suppose that each bridge graph has only one bridge edge, i.e., $|\mathcal{E}_{ij}^{S}|=1$.
The graph $\mathcal{G}(T)=\mathcal{G}^{S}\cup\mathcal{G}^{I}$ is SSC if each disjoint component $\mathcal{G}_{i}, i\in \{1,...,m\}$, is SSC.
\end{thm}

\begin{pf}
The \textit{if} condition can be proved by an induction of \textit{Corollary~\ref{corollary_1}}.
Let a state graph $\mathcal{G}^{S}$ be a pactus. 
Then, the pactus $\mathcal{G}^{S}$ can be induced by $m$-disjoint components $\mathcal{G}_{i}^{S}$ with $\mathcal{E}^{S}_{ij}$.
Also, each $\mathcal{G}_{i}^{S}$ is either a path or a cycle.
Suppose that each disjoint component $\mathcal{G}_{i}=\mathcal{G}_{i}^{S}\cup\mathcal{G}_{i}^{I}, i\in \{1,...,m\}$ satisfies
\textit{Proposition~\ref{proposition_path}} (path) or \textit{Theorem~\ref{theorem_cycle}} (cycle).
From each disjoint component point of view, it is clear that the set $\mathcal{N}(\alpha_{i})$$\setminus$$\alpha_{i}$ has at least one dedicated node 
for all $\alpha_{i}\subseteq\mathcal{V}_{i}^{S}, i\in \{1,...,m\}$.
Thus, since each SSC component $\mathcal{G}_{i}^{S}$ and $\mathcal{G}_{j}^{S}$ is connected by exactly one bridge edge $(k,l)\in\mathcal{E}_{ij}^{S}$
with $k\in\mathcal{V}^{S}_{i}$ and $l\in\mathcal{V}^{S}_{j}$ for $i\in \{1,...,m\}$ and $j\in\mathcal{N}_{\mathcal{G}_i^S}$.
By an induction of \textit{Corollary~\ref{corollary_1}}, 
the merged graph $\mathcal{G}(T)$ is SSC, 
i.e., the set $\mathcal{N}(\alpha)$$\setminus$$\alpha$ has at least one dedicated node
for all $\alpha\subseteq {\bigcup}^{m}_{i=1}\mathcal{V}_{i}^{S}$, 
which is equivalent to $\alpha\subseteq\mathcal{V}^{S}$.
\hfill $\square$
\end{pf}

The above \textit{Theorem~\ref{theorem_sym_pactus_sc_condition}} shows a sufficient condition for the strong structural controllability of pactus satisfying $|\mathcal{E}_{ij}^{S}|=1$, 
which is interpreted from the perspective of each component.
For example, let us consider the pactus depicted in Fig.~\ref{network_ex_th2}.
According to \textit{Proposition~\ref{proposition_path}}, $\mathcal{G}_{1}$ needs an external input node connected to node $1$ to be SSC, 
i.e., $\mathcal{V}_{1}^{I}=\{ u_{1}\}$. 
Since $\mathcal{G}_{2},\mathcal{G}_{3}$, and $\mathcal{G}_{4}$ are cycles, 
each component requires at least two properly located external input nodes to satisfy \textit{Theorem~\ref{theorem_cycle}},
i.e., $\mathcal{V}^{I}_{2}=\{ u_2,u_3\}$, $\mathcal{V}^{I}_{3}=\{ u_4,u_5\}$, $\mathcal{V}^{I}_{4}=\{ u_6,u_7\}$.
These results are shown in Fig.~\ref{network_ex_th2}. Hence, the graph $\mathcal{G}(T)=\mathcal{G}^{S}\cup\mathcal{G}^{I}$ requires seven external input nodes 
to satisfy \textit{Theorem~\ref{theorem_sym_pactus_sc_condition}},
i.e., $\mathcal{V}^{I}={\bigcup}^{4}_{i=1}\mathcal{V}_{i}^{I}=\{ u_{1},u_{2},u_{3},u_{4},u_{5},u_{6},u_{7}\}$. 
Note that the locations of the external input nodes are not unique.

\section{Strongly Structurally Controllable Graphs with Minimum External Input Nodes} \label{sec_topo_minimum}
In this section, we first present the condition for achieving strong structural controllability from the perspective of a single node.
Subsequently, a composition process of approximate polynomial complexity for MIP is provided.
For further analysis from a single node point of view, it is necessary to examine whether a state node is guaranteed at least one dedicated node.
For this reason, we introduce the notion of \textit{SSC node}, 
defined as a state node that satisfies the condition of strong structural controllability, a notion initially established in \cite{park2023KD}.

\begin{defn}\label{TC_state_node} 
\cite{park2023KD} (SSC node)
A state node $k\in\mathcal{V}_{i}^{S}$ is called a SSC node
if the set $\mathcal{N}(\alpha)$$\setminus$$\alpha$ has at least one dedicated node 
for all $\alpha\subseteq\mathcal{V}_{i}^{S}$ satisfying $k\in\alpha$.
The set of SSC nodes in $\mathcal{G}^{S}_{i}$ be symbolically written as $\mathcal{V}_{i}^{\text{\tiny SSC}}$.
\end{defn}
For convenience, we say that a state node $k\in\mathcal{V}^{S}$ has a dedicated node 
if the set $\mathcal{N}(\alpha)$$\setminus$$\alpha$ has at least one dedicated node 
for all $\alpha\subseteq\mathcal{V}^{S}$ satisfying $k\in\alpha$.
For example, consider the graph depicted in Fig.~\ref{network_ex_TCinput}(a).
In this case, the nodes $1,3\in\mathcal{V}^{S}$ are SSC nodes, 
which are guaranteed a dedicated node from the external input nodes $u_1$ and $u_2$, i.e., $1,3\in\mathcal{V}^{\text{\tiny SSC}}$.
In other words, if $\alpha$ contains at least one SSC node, the set $\mathcal{N}(\alpha)$$\setminus$$\alpha$ always has at least one dedicated node.
Based on the concept of an SSC node, 
the following corollary, which provides an interpretation of strong structural controllability from the perspective of each component, is directly derived from \cite{park2023KD}:

\begin{cor}\label{corollary_TCC}
Let us consider a pactus $\mathcal{G}^{S}=(\mathcal{V}^{S},\mathcal{E}^{S})$ induced by $m$-disjoint components 
$\mathcal{G}_{i}^{S}=({\mathcal{V}}_{i}^{S},{\mathcal{E}}_{i}^{S})$ with $\mathcal{E}_{ij}^{S}$ for $i\in \{1,...,m\}$ and $j\in\mathcal{N}_{\mathcal{G}_i^S}$.
Then, the graph $\mathcal{G}(T)=\mathcal{G}^{S}\cup\mathcal{G}^{I}$ is SSC
if and only if all state nodes in each component are SSC nodes, i.e., ${\bigcup}^{m}_{i=1}\mathcal{V}_{i}^{\text{\tiny SSC}}=\mathcal{V}^{S}$.
\end{cor}

In the context of MIP, a precise understanding of the role of input nodes is essential.
The following remark provides the significance of external input nodes, particularly from the perspective of dedicated nodes.

\begin{rem}\label{remark_meaning_input}
For a graph $\mathcal{G}(T)=(\mathcal{V},\mathcal{E})$, consider a state node $k\in\mathcal{V}^{S}$ that is connected to an external input node $u\in\mathcal{V}^{I}$ with an edge $(u,k)\in\mathcal{E}$.
Suppose a set $\alpha\subseteq\mathcal{V}^S$ includes the state node $k$, i.e., $k\in\alpha$. 
Given the assumption that an external input node is connected to exactly one state node, the node $u$ always meets the cardinality condition $|\mathcal{N}_{u}\cap\alpha|=1$.
Therefore, the external input node $u\in\mathcal{V}^{I}$ ensures the existence of a dedicated node for its connected state node $k\in\mathcal{V}^{S}$
\end{rem}

As shown in \textit{Remark~\ref{remark_meaning_input}}, 
the property of an external input node is to guarantee the existence of a dedicated node for any state node to which it is connected.
Surprisingly, this role of an external input node can also be fulfilled by a state node.
More specifically, if a specific structural condition is satisfied in a pactus, 
there exists a case that a state node $k\in\mathcal{V}_{i}^{S}$ in $\mathcal{G}_{i}$ guarantees 
the existence of a dedicated node of a state node $l\in\mathcal{V}_{j}^{S}$
in another component $\mathcal{G}_{j}$, which is adjacent to the node of $k$, i.e., $l\in\mathcal{N}_{k}$,
we call such state nodes \textit{component input nodes}.
With the above observation, the input nodes can be classified as the external input nodes and the component input nodes.


\begin{defn}\label{external_input_node} 
(External input node) A set of external input nodes in $\mathcal{G}_{i}$ is symbolically written as $\mathcal{V}_{i}^{IE}$. 
If a node $k\in\mathcal{V}_{i}^{I}$ guarantees a dedicated node of $l\in\mathcal{V}^{S}_{i}$ with a directed edge $(k,l)\in\mathcal{E}$,
the node $k$ is called an external input node of $\mathcal{G}_i$, i.e., $k\in\mathcal{V}^{IE}_{i}$.
\end{defn}

\begin{defn}\label{component_input_node} 
(Component input node) A set of component input nodes in $\mathcal{G}_{i}$ is symbolically written as $\mathcal{V}_{i}^{IC}$. 
Consider a graph $\mathcal{G}(T)=\mathcal{G}_{i}\cup\mathcal{G}_{ij}^{S}\cup\mathcal{G}_{j}$.
If a node $k\in\mathcal{V}_{i}^{S}$ guarantees a dedicated node of $l\in\mathcal{V}^{S}_{j}$ with an edge $(k,l)\in\mathcal{E}_{ij}^{S}$,
the node $k$ is called a component input node of $\mathcal{G}_j$, i.e., $k\in\mathcal{V}^{IC}_{j}$.
\end{defn}

\begin{figure}[]
\centering
\subfigure[]{
\begin{tikzpicture}[scale=0.4]
\node[] at (2.9,4.7) {\scriptsize$\mathcal{G}_1$};
\node[] at (7,6) {\scriptsize$\mathcal{G}_2$};

\node[red] at (4.7,4.2) {\scriptsize$\mathcal{G}_{12}^{S}$};

\node[place, black] (node1) at (1,5) [label=below:\scriptsize$1$] {};
\node[place, red] (node2) at (3,4) [label=below:\scriptsize$2$] {};
\node[place, red] (node3) at (5,3) [label=below:\scriptsize$3$] {};
\node[place, red] (node4) at (7,4) [label=below:\scriptsize$4$] {};
\node[place, black] (node5) at (5,5) [label=right:\scriptsize$5$] {};
\node[place, red] (node6) at (5,7) [label=above:\scriptsize$6$] {};

\node[place, black] (node7) at (7,8) [label=above:\scriptsize$7$] {};
\node[place, black] (node8) at (9,7) [label=above:\scriptsize$8$] {};
\node[place, black] (node9) at (9,5) [label=below:\scriptsize$9$] {};

\draw (node1) [line width=0.5pt] -- node [left] {} (node2);
\draw (node2) [line width=0.5pt] -- node [left] {} (node3);
\draw (node2) [red, dashed, line width=0.5pt] -- node [below] {} (node6);
\draw (node3) [red, dashed, line width=0.5pt] -- node [below] {} (node4);
\draw (node4) [line width=0.5pt] -- node [left] {} (node5);
\draw (node5) [line width=0.5pt] -- node [left] {} (node6);
\draw (node6) [line width=0.5pt] -- node [left] {} (node7);
\draw (node7) [line width=0.5pt] -- node [left] {} (node8);
\draw (node8) [line width=0.5pt] -- node [left] {} (node9);
\draw (node4) [line width=0.5pt] -- node [left] {} (node9);

\node[place, circle] (node17) at (-0.5,4) [label=below:\scriptsize$u_{1}$] {}; 
\node[place, circle] (node18) at (3.5,2) [label=below:\scriptsize$u_{2}$] {}; 
\draw (node17) [-latex, line width=0.5pt] -- node [right] {} (node1);
\draw (node18) [-latex, line width=0.5pt] -- node [right] {} (node3);
\end{tikzpicture}
}
\subfigure[]{
\begin{tikzpicture}[scale=0.43]
\node[] at (7,6) {\scriptsize$\mathcal{G}_2$};

\node[place, black] (node4) at (7,4) [label=below:\scriptsize$4$] {};
\node[place, black] (node5) at (5,5) [label=below:\scriptsize$5$] {};
\node[place, black] (node6) at (5,7) [label=above:\scriptsize$6$] {};

\node[place, black] (node7) at (7,8) [label=above:\scriptsize$7$] {};
\node[place, black] (node8) at (9,7) [label=above:\scriptsize$8$] {};
\node[place, black] (node9) at (9,5) [label=below:\scriptsize$9$] {};

\draw (node4) [line width=0.5pt] -- node [left] {} (node5);
\draw (node5) [line width=0.5pt] -- node [left] {} (node6);
\draw (node6) [line width=0.5pt] -- node [left] {} (node7);
\draw (node7) [line width=0.5pt] -- node [left] {} (node8);
\draw (node8) [line width=0.5pt] -- node [left] {} (node9);
\draw (node4) [line width=0.5pt] -- node [left] {} (node9);

\node[place, circle] (node17) at (5,3) [label=below:\scriptsize$3$] {}; 
\node[place, circle] (node18) at (3,4) [label=below:\scriptsize$2$] {}; 
\draw (node4) [line width=0.5pt] -- node [left] {} (node17);
\draw (node6) [line width=0.5pt] -- node [left] {} (node18);
\end{tikzpicture}
}
\caption{Example of the SSC nodes and the component input nodes:
(a) A graph $\mathcal{G}(T)=\mathcal{G}_{1}\cup\mathcal{G}_{12}\cup\mathcal{G}_{2}$ with $|\mathcal{V}^{IE}|=2$. 
(b) $\mathcal{G}_{2}$ with $\mathcal{V}^{IC}_{2}=\{ 2,3 \}$.
}
\label{network_ex_TCinput}
\end{figure}

From the above definitions, the set of input nodes $\mathcal{V}^{I}$ is re-defined as 
a union of the set of component input nodes and the set of external input nodes, 
i.e., $\mathcal{V}^{I}=\mathcal{V}^{IE}\cup\mathcal{V}^{IC}$ satisfying $\mathcal{V}^{IE}\cap\mathcal{V}^{IC}=\emptyset$.
For a pactus given by $\mathcal{G}^{S}={\bigcup}^{m}_{i=1}{\bigcup}_{j\in\mathcal{N}_{\mathcal{G}_i^S}}(\mathcal{G}_{i}^{S}    \cup    \mathcal{G}_{ij}^{S})$,
The following theorem provides the necessary and sufficient condition for the component input nodes in a graph $\mathcal{G}(T)=\mathcal{G}^S\cup\mathcal{G}^I$.

\begin{thm}\label{theorem_TCC}
Consider a graph $\mathcal{G}(T)=\mathcal{G}^S\cup\mathcal{G}^I$ with a pactus $\mathcal{G}^{S}$.
Then, the state nodes in $\mathcal{V}_{i}^{S}\cap \mathcal{V}_{ij}^{S}$, are the component input nodes of $\mathcal{G}_{j}$, 
i.e., $\mathcal{V}_{j}^{IC}=\mathcal{V}_{i}^{S}\cap\mathcal{V}_{ij}^{S}$,
if and only if ${\bigcup}_{j\in\mathcal{N}_{\mathcal{G}_i^S}}(\mathcal{G}_{i}\cup\mathcal{G}_{ij}^{S})$ is SSC for $i\in\{1,...,m\}$.
\end{thm}

\begin{pf}
Let the subgraph ${\bigcup}_{j\in\mathcal{N}_{\mathcal{G}_i^S}}(\mathcal{G}_{i}\cup\mathcal{G}_{ij}^{S})$ be denoted as 
$\bar{\mathcal{G}}_i=(\bar{\mathcal{V}}_i,\bar{\mathcal{E}}_i)$ for $i\in\{1,...,m\}$.
For \textit{if} condition, suppose that $\bar{\mathcal{G}}_i$ is SSC for each $i\in \{1,...,m\}$.
According to \textit{Definition~\ref{TC_state_node}} and \textit{Corollary~\ref{corollary_TCC}}, it follows that all state nodes in $\bar{\mathcal{V}}_{i}$ are SSC nodes.
Obviously, any bridge nodes $k,l\in\mathcal{V}_{ij}^{S}$ satisfying $k\in\mathcal{N}_{l}$ are also SSC nodes.
This ensures that the node $k\in\mathcal{V}^{S}_{i}$ guarantees the existence of a dedicated node of the node $l\in\mathcal{V}^{S}_{j}$.
Hence, by \textit{Definition~\ref{component_input_node}}, the bridge nodes in $\mathcal{V}_{i}^{S}\cap\mathcal{V}_{ij}^{S}$ are the component input nodes of $\mathcal{G}_{j}$,
i.e., $\mathcal{V}_{j}^{IC}=\mathcal{V}_{i}^{S}\cap\mathcal{V}_{ij}^{S}$.

For the \textit{only if} condition, suppose $\bar{\mathcal{G}}_i$ is not SSC for some $i\in \{1,...,m\}$. 
Note that if $\bar{\mathcal{G}}_i$ is SSC, then $\mathcal{G}_{i}$ must also be SSC, but the converse is not necessarily true. 
Let us consider the case where $\mathcal{G}_{i}$ is SSC. 
Under this assumption, the bridge nodes in $\mathcal{V}_{i}^{S}\cap\mathcal{V}_{ij}^{S}$ are SSC nodes. 
However, the assumption that $\bar{\mathcal{G}}_i$ is not SSC implies that there exists at least one bridge node in $\mathcal{V}_{j}^{S}\cap\mathcal{V}_{ij}^{S}$ that is not an SSC node. 
This results in a case where there is no dedicated node for some subsets $\alpha\subseteq\mathcal{V}^{S}$ including a node in $\mathcal{V}_{j}^{S}\cap\mathcal{V}_{ij}^{S}$.
\hfill $\square$
\end{pf}

Note that if \textit{Theorem~\ref{theorem_TCC}} is satisfied, the set of component input nodes is defined as 
$\mathcal{V}^{IC}_{j}=\mathcal{V}^{S}_{i}\cap\mathcal{V}^{S}_{ij}$ for $i\in \{1,...,m\}$ and $j\in\mathcal{N}_{\mathcal{G}_i^S}$.
This implies that the component input nodes of $\mathcal{G}_j$ can be determined only by its adjacent component $\mathcal{G}_i,i\in\mathcal{N}_{\mathcal{G}_j^S}$.
For example, in the graph $\mathcal{G}(T)$ depicted in Fig.~\ref{network_ex_TCinput}(a), 
consider the SSC subgraph $\mathcal{G}_{1}\cup\mathcal{G}_{12}^{S}$ satisfying \textit{Theorem~\ref{theorem_tree}}.
Then, according to \textit{Theorem~\ref{theorem_TCC}}, the nodes $2,3\in\mathcal{V}_{1}^{S}\cap\mathcal{V}_{12}^{S}$ are 
the component input nodes of $\mathcal{G}_{2}$, i.e., $\mathcal{V}_{2}^{IC}=\{2,3\}$.
From the viewpoint of $\mathcal{G}_{2}$, the graph depicted in Fig.~\ref{network_ex_TCinput}(a) can be expressed as shown in Fig.~\ref{network_ex_TCinput}(b).
Note that the property of an external input node and component input node is equivalent to the existence of a dedicated node.
Now, let us suppose that the graph $\mathcal{G}_{2}$ in Fig.~\ref{network_ex_TCinput}(b) satisfies \textit{Theorem~\ref{theorem_cycle}}
by an additional properly located external input node connected to one of the nodes $5,7$, and $9$.
Then, all state nodes in $\mathcal{G}(T)$ become SSC nodes, it follows from \textit{Corollary~\ref{corollary_TCC}} that $\mathcal{G}(T)$ is SSC.
In this manner, based on the concept of component input node,
SSC graphs with the minimum number of external input nodes can be designed by adding additional proper external input nodes.

\begin{rem}\label{remark_IC}
While most existing literature on the MIP for structured networks primarily focuses on external input nodes, our paper introduces a novel concept of component input nodes.
Despite being state nodes, these component input nodes function identically to external input nodes. 
Based on the notions of dedicated and sharing nodes, we demonstrate that these nodes indeed assure the presence of dedicated nodes, thereby playing a role equivalent to external input nodes. 
It is important to note that this finding does not invalidate previous studies concentrated on external input nodes for the MIP. 
Instead, it implies that in certain graph structures, control paths emanating from external input nodes may inherently incorporate the component input nodes we define. 
This insight provides a fresh perspective on the topological formation of control paths, offering significant contributions and broadening existing methodologies in the field.
\end{rem}

For a graph $\mathcal{G}(T)=\mathcal{G}^{S}\cup\mathcal{G}^{I}$, the minimum number of external input nodes in $\mathcal{V}^{I}$ 
for the strong structural controllability is symbolically written as $\min{|\mathcal{V}^{IE}|}$.
Note that it follows from \textit{Proposition~\ref{proposition_path}}, \textit{Theorem~\ref{theorem_tree}}, and \textit{Theorem~\ref{theorem_cycle}} that 
$\min{|\mathcal{V}^{I}|}$ of path, tree, and cycle are $1$, $m$, and $2$, respectively.
For a general type of pactus, here we propose a composition algorithm of polynomial complexity
that can uniquely determine the minimum number of external input nodes while maintaining strong structural controllability.

The algorithm starts with a \textit{decomposition process}, which involves breaking down a pactus into its basic components.
This process sequentially classifies the structure of the given pactus, 
represented as $\mathcal{G}^{S}={\bigcup}^{m}_{i=1}{\bigcup}_{j\in\mathcal{N}_{\mathcal{G}_i^S}}(\mathcal{G}_{i}^{S}    \cup    \mathcal{G}_{ij}^{S})$, into specific types in the decomposition process.
In this \textit{decomposition process}, the graph $\bar{\mathcal{G}}_{i}^{S}$ denotes a union of state graphs of the $i$-th component and its bridge graph,
i.e., $\bar{\mathcal{G}}_{i}^{S} := \mathcal{G}_{i}^{S}\cup\mathcal{G}_{i(i+1)}^{S}$, where $i+1\in\mathcal{N}_{\mathcal{G}_i^S}$.
We classify graph types of  $\bar{\mathcal{G}}_{i}^{S}$ as \textit{path-type}, \textit{tree-type}, and \textit{cycle-type}.
The graph $\bar{\mathcal{G}}_{i}^{S}$ is \textit{path-type} and \textit{tree-type} if $\bar{\mathcal{G}}_{i}^{S}$ is a path and a tree graph, respectively.
Furthermore, $\bar{\mathcal{G}}_{i}^{S}$ is \textit{cycle-type} if $\bar{\mathcal{G}}_{i}^{S}$ contains a cycle graph.
The algorithm continues with a \textit{composition process}, which involves methodically adding the minimum number of external input nodes to each component. 
This process is conducted step-by-step, following the guidelines established in \textit{Proposition~\ref{proposition_path}}, \textit{Theorem~\ref{theorem_tree}}, and \textit{Theorem~\ref{theorem_cycle}}.
Accordingly, based on \textit{Theorem~\ref{theorem_TCC}}, the set of component input nodes for each stage is determined as 
$\mathcal{V}_{i+1}^{IC} :=\mathcal{V}_{i}^{S}\cap\mathcal{V}^{S}_{i(i+1)}$, where $i+1\in\mathcal{N}_{\mathcal{G}_i^S}$.
Utilizing the above principles, the graph composition algorithm for the general type of a pactus is outlined in \textit{Algorithm~\ref{graph_composition}}.

\begin{algorithm}[h] 
\caption{}\label{graph_composition}
\begin{algorithmic}[1]
\BState \textbf{Decomposition}:
\State a pactus 
$\mathcal{G}^{S}:={\bigcup}^{m}_{i=1}{\bigcup}_{j\in\mathcal{N}_{\mathcal{G}_i^S}}(\mathcal{G}_{i}^{S}    \cup    \mathcal{G}_{ij}^{S})$, $\mathcal{V}^{IE}=\emptyset$ 
\State $i=0$
 \For{$i=i+1$} 
     \If { $i < m$}
    \State $\bar{\mathcal{G}}_{i}^{S} := \mathcal{G}_{i}^{S}\cup\mathcal{G}_{i(i+1)}^{S}$ where $i+1\in\mathcal{N}_{\mathcal{G}_i^S}$
    \ElsIf { $i = m$}
    \State $\bar{\mathcal{G}}_{i} := \mathcal{G}_{i}^{S}$
    \Else 
\State \textit{goto} \textbf{composition}
    \EndIf
  \EndFor   \Comment{End for $i$}\newline

\BState \textbf{composition}:
\State $i=0$
 \For{$i=i+1$} 

\If { $\bar{\mathcal{G}}_{i}^{S}$ is not SSC} 
  \Switch{graph type of $\bar{\mathcal{G}}_{i}^{S}$}
\Case{\textbf{path-type}}
\State add $\min{|\mathcal{V}^{IE}_{i}|}$ to satisfy \textit{Proposition~\ref{proposition_path}}
\EndCase

\Case{\textbf{tree-type}}
\State add $\min{|\mathcal{V}^{IE}_{i}|}$ to satisfy \textit{Theorem~\ref{theorem_tree}}
\EndCase

\Case{\textbf{cycle-type}}
\State add $\min{|\mathcal{V}^{IE}_{i}|}$ to satisfy \textit{Theorem~\ref{theorem_cycle}}
\EndCase
  \EndSwitch
\EndIf

\If { $i < m$}
    \State update $\mathcal{V}_{i+1}^{IC} :=\mathcal{V}_{i}^{S}\cap\mathcal{V}^{S}_{i(i+1)}$
    \State $\bar{\mathcal{G}}_{i}:=\bar{\mathcal{G}}_{i}^{S}\cup\bar{\mathcal{G}}_{i}^{I}$
    \Else 
    \State 
     \textit{goto} \textit{end}
    \EndIf
  \EndFor   \Comment{End for $i$}
\BState \emph{end}:
\State \textit{output} $\mathcal{G}(T)^{out}=\bigcup^{m}_{i=1}\bar{\mathcal{G}}_{i}$  with $min|\mathcal{V}^{IE}|$
\end{algorithmic}
\end{algorithm}

\begin{figure}[]
\centering
\subfigure[Initial : A pactus $\mathcal{G}^{S}$]{
\begin{tikzpicture}[scale=0.45]
\node[] at (2.9,4.7) {\scriptsize$\mathcal{G}_1$};
\node[] at (7,6) {\scriptsize$\mathcal{G}_2$};
\node[] at (10.4,3) {\scriptsize$\mathcal{G}_3$};
\node[] at (14,4.5) {\scriptsize$\mathcal{G}_4$};

\node[] at (4.7,4.2) {\scriptsize$\mathcal{G}_{12}$};
\node[] at (9,4) {\scriptsize$\mathcal{G}_{23}$};
\node[] at (12.2,4) {\scriptsize$\mathcal{G}_{34}$};

\node[place, black] (node1) at (1,5) [label=below:\scriptsize$1$] {};
\node[place, black] (node2) at (3,4) [label=below:\scriptsize$2$] {};
\node[place, black] (node3) at (5,3) [label=below:\scriptsize$3$] {};
\node[place, black] (node4) at (7,4) [label=below:\scriptsize$4$] {};
\node[place, black] (node5) at (5,5) [label=right:\scriptsize$5$] {};
\node[place, black] (node6) at (5,7) [label=above:\scriptsize$6$] {};
\node[place, black] (node7) at (7,8) [label=above:\scriptsize$7$] {};
\node[place, black] (node8) at (9,7) [label=above:\scriptsize$8$] {};
\node[place, black] (node9) at (9,5) [label=left:\scriptsize$9$] {};
\node[place, black] (node10) at (9,3) [label=below:\scriptsize$10$] {};
\node[place, black] (node11) at (11,2) [label=below:\scriptsize$11$] {};
\node[place, black] (node12) at (11,4) [label=above:\scriptsize$12$] {};
\node[place, black] (node13) at (13,5) [label=above:\scriptsize$13$] {};
\node[place, black] (node14) at (15,6) [label=above:\scriptsize$14$] {};
\node[place, black] (node15) at (15,4) [label=below:\scriptsize$15$] {};
\node[place, black] (node16) at (13,3) [label=below:\scriptsize$16$] {};

\draw (node1) [line width=0.5pt] -- node [left] {} (node2);
\draw (node2) [line width=0.5pt] -- node [left] {} (node3);
\draw (node2) [line width=0.5pt] -- node [below] {} (node6);
\draw (node3) [,line width=0.5pt] -- node [below] {} (node4);
\draw (node4) [line width=0.5pt] -- node [left] {} (node5);
\draw (node5) [line width=0.5pt] -- node [left] {} (node6);
\draw (node6) [line width=0.5pt] -- node [left] {} (node7);
\draw (node7) [line width=0.5pt] -- node [left] {} (node8);
\draw (node8) [line width=0.5pt] -- node [left] {} (node9);
\draw (node4) [line width=0.5pt] -- node [left] {} (node9);
\draw (node4) [line width=0.5pt] -- node [below] {} (node10);
\draw (node10) [line width=0.5pt] -- node [left] {} (node11);
\draw (node10) [line width=0.5pt] -- node [left] {} (node12);
\draw (node11) [line width=0.5pt] -- node [left] {} (node12);
\draw (node9) [line width=0.5pt] -- node [below] {} (node12);
\draw (node12) [line width=0.5pt] -- node [below] {} (node13);
\draw (node13) [line width=0.5pt] -- node [left] {} (node14);
\draw (node14) [line width=0.5pt] -- node [left] {} (node15);
\draw (node15) [line width=0.5pt] -- node [left] {} (node16);
\draw (node13) [line width=0.5pt] -- node [left] {} (node16);
\end{tikzpicture}
}
\subfigure[Step 1 : $\bar{\mathcal{G}}_{1}=\mathcal{G}_{1}\cup\mathcal{G}_{12}$]{
\begin{tikzpicture}[scale=0.45]

\node[] at (2.9,4.7) {\scriptsize$\mathcal{G}_1$};
\node[] at (7,6) {\scriptsize$\mathcal{G}_2$};
\node[] at (10.4,3) {\scriptsize$\mathcal{G}_3$};
\node[] at (14,4.5) {\scriptsize$\mathcal{G}_4$};

\node[] at (4.7,4.2) {\scriptsize$\mathcal{G}_{12}$};
\node[] at (9,4) {\scriptsize$\mathcal{G}_{23}$};
\node[] at (12.2,4) {\scriptsize$\mathcal{G}_{34}$};

\node[place, blue] (node1) at (1,5) [label=below:\scriptsize$1$] {};
\node[place, blue] (node2) at (3,4) [label=below:\scriptsize$2$] {};
\node[place, blue] (node3) at (5,3) [label=below:\scriptsize$3$] {};
\node[place, blue] (node4) at (7,4) [label=below:\scriptsize$4$] {};
\node[place, black] (node5) at (5,5) [label=right:\scriptsize$5$] {};
\node[place, blue] (node6) at (5,7) [label=above:\scriptsize$6$] {};
\node[place, black] (node7) at (7,8) [label=above:\scriptsize$7$] {};
\node[place, black] (node8) at (9,7) [label=above:\scriptsize$8$] {};
\node[place, black] (node9) at (9,5) [label=left:\scriptsize$9$] {};
\node[place, black] (node10) at (9,3) [label=below:\scriptsize$10$] {};
\node[place, black] (node11) at (11,2) [label=below:\scriptsize$11$] {};
\node[place, black] (node12) at (11,4) [label=above:\scriptsize$12$] {};
\node[place, black] (node13) at (13,5) [label=above:\scriptsize$13$] {};
\node[place, black] (node14) at (15,6) [label=above:\scriptsize$14$] {};
\node[place, black] (node15) at (15,4) [label=below:\scriptsize$15$] {};
\node[place, black] (node16) at (13,3) [label=below:\scriptsize$16$] {};

\draw (node1) [line width=0.5pt] -- node [left] {} (node2);
\draw (node2) [line width=0.5pt] -- node [left] {} (node3);
\draw (node2) [line width=0.5pt] -- node [below] {} (node6);
\draw (node3) [line width=0.5pt] -- node [below] {} (node4);
\draw (node4) [line width=0.5pt] -- node [left] {} (node5);
\draw (node5) [line width=0.5pt] -- node [left] {} (node6);
\draw (node6) [line width=0.5pt] -- node [left] {} (node7);
\draw (node7) [line width=0.5pt] -- node [left] {} (node8);
\draw (node8) [line width=0.5pt] -- node [left] {} (node9);
\draw (node4) [line width=0.5pt] -- node [left] {} (node9);
\draw (node4) [line width=0.5pt] -- node [below] {} (node10);
\draw (node10) [line width=0.5pt] -- node [left] {} (node11);
\draw (node10) [line width=0.5pt] -- node [left] {} (node12);
\draw (node11) [line width=0.5pt] -- node [left] {} (node12);
\draw (node9) [line width=0.5pt] -- node [below] {} (node12);
\draw (node12) [line width=0.5pt] -- node [below] {} (node13);
\draw (node13) [line width=0.5pt] -- node [left] {} (node14);
\draw (node14) [line width=0.5pt] -- node [left] {} (node15);
\draw (node15) [line width=0.5pt] -- node [left] {} (node16);
\draw (node13) [line width=0.5pt] -- node [left] {} (node16);

\node[place, circle] (node17) at (-0.5,4) [label=below:\scriptsize$u_{1}$] {}; 
\node[place, circle] (node18) at (3.5,2) [label=below:\scriptsize$u_{2}$] {}; 
\draw (node17) [-latex, line width=0.5pt] -- node [right] {} (node1);
\draw (node18) [-latex, line width=0.5pt] -- node [right] {} (node3);
\end{tikzpicture}
}
\subfigure[Step 2 : $\bar{\mathcal{G}}_{2}=\mathcal{G}_{2}\cup\mathcal{G}_{23}$]{
\begin{tikzpicture}[scale=0.45]

\node[] at (7,6) {\scriptsize$\mathcal{G}_2$};
\node[] at (10.4,3) {\scriptsize$\mathcal{G}_3$};
\node[] at (14,4.5) {\scriptsize$\mathcal{G}_4$};

\node[] at (4.7,4.2) {\scriptsize$\mathcal{G}_{12}$};
\node[] at (9,4) {\scriptsize$\mathcal{G}_{23}$};
\node[] at (12.2,4) {\scriptsize$\mathcal{G}_{34}$};

\node[place, circle] (node2) at (3,4) [label=below:\scriptsize$2$] {};
\node[place, circle] (node3) at (5,3) [label=below:\scriptsize$3$] {};
\node[place, blue] (node4) at (7,4) [label=below:\scriptsize$4$] {};
\node[place, blue] (node5) at (5,5) [label=right:\scriptsize$5$] {};
\node[place, blue] (node6) at (5,7) [label=above:\scriptsize$6$] {};
\node[place, blue] (node7) at (7,8) [label=above:\scriptsize$7$] {};
\node[place, blue] (node8) at (9,7) [label=above:\scriptsize$8$] {};
\node[place, blue] (node9) at (9,5) [label=left:\scriptsize$9$] {};
\node[place, blue] (node10) at (9,3) [label=below:\scriptsize$10$] {};
\node[place, black] (node11) at (11,2) [label=below:\scriptsize$11$] {};
\node[place, blue] (node12) at (11,4) [label=above:\scriptsize$12$] {};
\node[place, black] (node13) at (13,5) [label=above:\scriptsize$13$] {};
\node[place, black] (node14) at (15,6) [label=above:\scriptsize$14$] {};
\node[place, black] (node15) at (15,4) [label=below:\scriptsize$15$] {};
\node[place, black] (node16) at (13,3) [label=below:\scriptsize$16$] {};

\draw (node2) [line width=0.5pt] -- node [below] {} (node6);
\draw (node3) [line width=0.5pt] -- node [below] {} (node4);
\draw (node4) [line width=0.5pt] -- node [left] {} (node5);
\draw (node5) [line width=0.5pt] -- node [left] {} (node6);
\draw (node6) [line width=0.5pt] -- node [left] {} (node7);
\draw (node7) [line width=0.5pt] -- node [left] {} (node8);
\draw (node8) [line width=0.5pt] -- node [left] {} (node9);
\draw (node4) [line width=0.5pt] -- node [left] {} (node9);
\draw (node4) [line width=0.5pt] -- node [below] {} (node10);
\draw (node10) [line width=0.5pt] -- node [left] {} (node11);
\draw (node10) [line width=0.5pt] -- node [left] {} (node12);
\draw (node11) [line width=0.5pt] -- node [left] {} (node12);
\draw (node9) [line width=0.5pt] -- node [below] {} (node12);
\draw (node12) [line width=0.5pt] -- node [below] {} (node13);
\draw (node13) [line width=0.5pt] -- node [left] {} (node14);
\draw (node14) [line width=0.5pt] -- node [left] {} (node15);
\draw (node15) [line width=0.5pt] -- node [left] {} (node16);
\draw (node13) [line width=0.5pt] -- node [left] {} (node16);

\node[place, circle] (node20) at (5.5,9) [label=below:\scriptsize$u_{3}$] {}; 
\draw (node20) [-latex, line width=0.5pt] -- node [right] {} (node7);
\end{tikzpicture}
}

\subfigure[Step 3 : $\bar{\mathcal{G}}_{3}=\mathcal{G}_{3}\cup\mathcal{G}_{34}$]{
\begin{tikzpicture}[scale=0.45]

\node[] at (10.4,3) {\scriptsize$\mathcal{G}_3$};
\node[] at (14,4.5) {\scriptsize$\mathcal{G}_4$};

\node[] at (9,4) {\scriptsize$\mathcal{G}_{23}$};
\node[] at (12.2,4) {\scriptsize$\mathcal{G}_{34}$};

\node[place, circle] (node4) at (7,4) [label=above:\scriptsize$4$] {};
\node[place, circle] (node9) at (9,5) [label=above:\scriptsize$9$] {};
\node[place, blue] (node10) at (9,3) [label=below:\scriptsize$10$] {};
\node[place, blue] (node11) at (11,2) [label=below:\scriptsize$11$] {};
\node[place, blue] (node12) at (11,4) [label=above:\scriptsize$12$] {};
\node[place, blue] (node13) at (13,5) [label=above:\scriptsize$13$] {};
\node[place, black] (node14) at (15,6) [label=above:\scriptsize$14$] {};
\node[place, black] (node15) at (15,4) [label=below:\scriptsize$15$] {};
\node[place, black] (node16) at (13,3) [label=below:\scriptsize$16$] {};

\draw (node4) [line width=0.5pt] -- node [below] {} (node10);
\draw (node10) [line width=0.5pt] -- node [left] {} (node11);
\draw (node10) [line width=0.5pt] -- node [left] {} (node12);
\draw (node11) [line width=0.5pt] -- node [left] {} (node12);
\draw (node9) [line width=0.5pt] -- node [below] {} (node12);
\draw (node12) [line width=0.5pt] -- node [below] {} (node13);
\draw (node13) [line width=0.5pt] -- node [left] {} (node14);
\draw (node14) [line width=0.5pt] -- node [left] {} (node15);
\draw (node15) [line width=0.5pt] -- node [left] {} (node16);
\draw (node13) [line width=0.5pt] -- node [left] {} (node16);

\end{tikzpicture}

}
\subfigure[Step 4 : $\bar{\mathcal{G}}_{4}=\mathcal{G}_{4}$]{
\begin{tikzpicture}[scale=0.45]

\node[] at (14,4.5) {\scriptsize$\mathcal{G}_4$};
\node[] at (16,4.5) {\scriptsize$ $};

\node[] at (12.2,4) {\scriptsize$\mathcal{G}_{34}$};

\node[place, circle] (node12) at (11,4) [label=above:\scriptsize$12$] {};
\node[place, blue] (node13) at (13,5) [label=above:\scriptsize$13$] {};
\node[place, blue] (node14) at (15,6) [label=above:\scriptsize$14$] {};
\node[place, blue] (node15) at (15,4) [label=below:\scriptsize$15$] {};
\node[place, blue] (node16) at (13,3) [label=below:\scriptsize$16$] {};

\draw (node12) [line width=0.5pt] -- node [below] {} (node13);
\draw (node13) [line width=0.5pt] -- node [left] {} (node14);
\draw (node14) [line width=0.5pt] -- node [left] {} (node15);
\draw (node15) [line width=0.5pt] -- node [left] {} (node16);
\draw (node13) [line width=0.5pt] -- node [left] {} (node16);

\node[place, circle] (node19) at (13.5,7) [label=below:\scriptsize$u_{4}$] {}; 
\draw (node19) [-latex, line width=0.5pt] -- node [right] {} (node14);
\end{tikzpicture}
}

\subfigure[Final : $\mathcal{G}(T)^{out}=\bar{\mathcal{G}}_{1}\cup\bar{\mathcal{G}}_{2}\cup\bar{\mathcal{G}}_{3}\cup\bar{\mathcal{G}}_{4}$]{
\begin{tikzpicture}[scale=0.45]

\node[] at (2.9,4.7) {\scriptsize$\mathcal{G}_1$};
\node[] at (7,6) {\scriptsize$\mathcal{G}_2$};
\node[] at (10.4,3) {\scriptsize$\mathcal{G}_3$};
\node[] at (14,4.5) {\scriptsize$\mathcal{G}_4$};

\node[] at (4.7,4.2) {\scriptsize$\mathcal{G}_{12}$};
\node[] at (9,4) {\scriptsize$\mathcal{G}_{23}$};
\node[] at (12.2,4) {\scriptsize$\mathcal{G}_{34}$};

\node[place, blue] (node1) at (1,5) [label=below:\scriptsize$1$] {};
\node[place, blue] (node2) at (3,4) [label=below:\scriptsize$2$] {};
\node[place, blue] (node3) at (5,3) [label=below:\scriptsize$3$] {};
\node[place, blue] (node4) at (7,4) [label=below:\scriptsize$4$] {};
\node[place, blue] (node5) at (5,5) [label=right:\scriptsize$5$] {};
\node[place, blue] (node6) at (5,7) [label=above:\scriptsize$6$] {};
\node[place, blue] (node7) at (7,8) [label=above:\scriptsize$7$] {};
\node[place, blue] (node8) at (9,7) [label=above:\scriptsize$8$] {};
\node[place, blue] (node9) at (9,5) [label=left:\scriptsize$9$] {};
\node[place, blue] (node10) at (9,3) [label=below:\scriptsize$10$] {};
\node[place, blue] (node11) at (11,2) [label=below:\scriptsize$11$] {};
\node[place, blue] (node12) at (11,4) [label=above:\scriptsize$12$] {};
\node[place, blue] (node13) at (13,5) [label=above:\scriptsize$13$] {};
\node[place, blue] (node14) at (15,6) [label=above:\scriptsize$14$] {};
\node[place, blue] (node15) at (15,4) [label=below:\scriptsize$15$] {};
\node[place, blue] (node16) at (13,3) [label=below:\scriptsize$16$] {};

\draw (node1) [line width=0.5pt] -- node [left] {} (node2);
\draw (node2) [line width=0.5pt] -- node [left] {} (node3);
\draw (node2) [line width=0.5pt] -- node [below] {} (node6);
\draw (node3) [,line width=0.5pt] -- node [below] {} (node4);
\draw (node4) [line width=0.5pt] -- node [left] {} (node5);
\draw (node5) [line width=0.5pt] -- node [left] {} (node6);
\draw (node6) [line width=0.5pt] -- node [left] {} (node7);
\draw (node7) [line width=0.5pt] -- node [left] {} (node8);
\draw (node8) [line width=0.5pt] -- node [left] {} (node9);
\draw (node4) [line width=0.5pt] -- node [left] {} (node9);
\draw (node4) [line width=0.5pt] -- node [below] {} (node10);
\draw (node10) [line width=0.5pt] -- node [left] {} (node11);
\draw (node10) [line width=0.5pt] -- node [left] {} (node12);
\draw (node11) [line width=0.5pt] -- node [left] {} (node12);
\draw (node9) [line width=0.5pt] -- node [below] {} (node12);
\draw (node12) [line width=0.5pt] -- node [below] {} (node13);
\draw (node13) [line width=0.5pt] -- node [left] {} (node14);
\draw (node14) [line width=0.5pt] -- node [left] {} (node15);
\draw (node15) [line width=0.5pt] -- node [left] {} (node16);
\draw (node13) [line width=0.5pt] -- node [left] {} (node16);
 
\node[place, circle] (node17) at (-0.5,4) [label=below:\scriptsize$u_{1}$] {}; 
\node[place, circle] (node18) at (3.5,2) [label=below:\scriptsize$u_{2}$] {}; 
\node[place, circle] (node20) at (5.5,9) [label=below:\scriptsize$u_{3}$] {}; 
\node[place, circle] (node19) at (13.5,7) [label=below:\scriptsize$u_{4}$] {}; 
\draw (node17) [-latex, line width=0.5pt] -- node [right] {} (node1); 
\draw (node18) [-latex, line width=0.5pt] -- node [right] {} (node3);
\draw (node20) [-latex, line width=0.5pt] -- node [right] {} (node7);
\draw (node19) [-latex, line width=0.5pt] -- node [right] {} (node14);
\end{tikzpicture}
}
\caption{Topological example of \textit{Algorithm~\ref{graph_composition}},
where the blue-marked nodes represent the SSC nodes.}
\label{network_ex_step}
\end{figure}
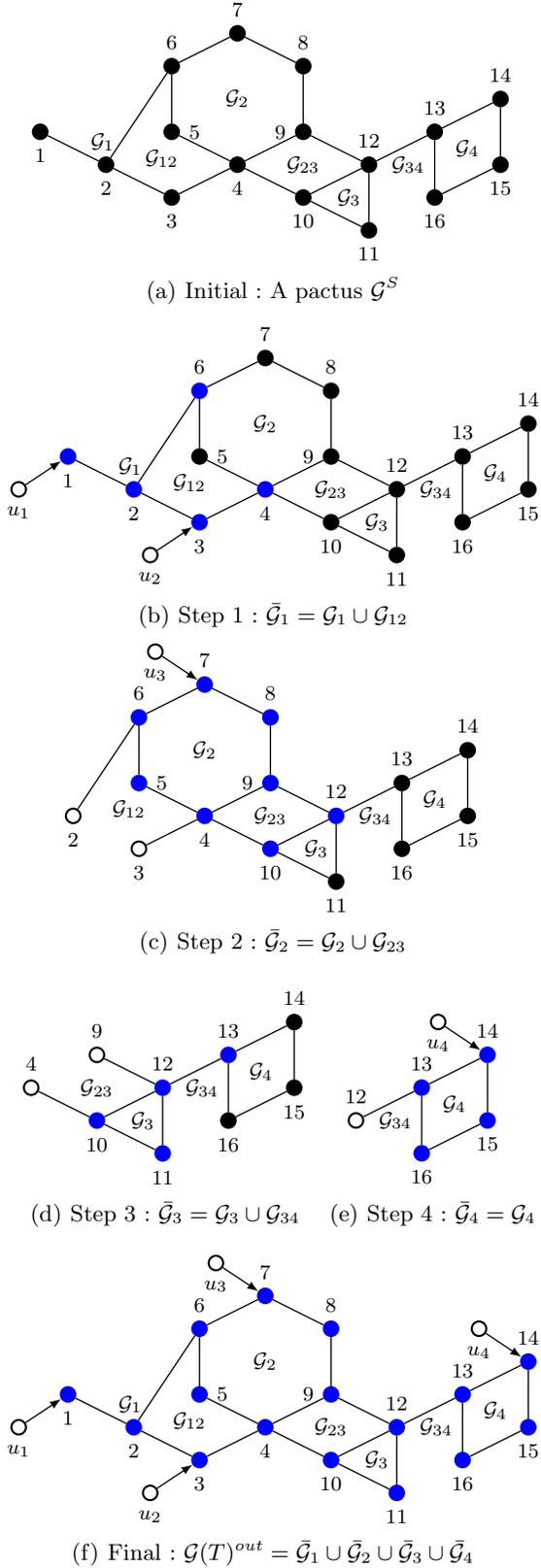

As a topological example of \textit{Algorithm~\ref{graph_composition}}, 
let us consider the pactus $\mathcal{G}^{S}=(\mathcal{V}^{S},\mathcal{E}^{S})$ depicted in Fig.~\ref{network_ex_step}(a) consisting of four disjoint components with bridge edges.
In \textit{decomposition process}, the pactus $\mathcal{G}^{S}$ is decomposed into disjoint components and bridge graphs.
In \textit{composition process}, for each $\bar{\mathcal{G}}_{i}=\bar{\mathcal{G}}_{i}^{S}\cup\bar{\mathcal{G}}_{i}^{I}, i\in\{1,2,3,4\}$, 
additional external input nodes are added at proper locations to guarantee the existence of dedicated nodes for all $\alpha\subseteq\bar{\mathcal{V}}_{i}^{S}$,
i.e, the SSC nodes in $\mathcal{G}$ will expand sequentially. 
As a result, the output graph $\mathcal{G}^{out}$ in \textit{Algorithm~\ref{graph_composition}} satisfies \textit{Corollary~\ref{corollary_TCC}}, 
which is equivalent to satisfying \textit{Theorem~\ref{theorem_Tsatsomeros}}.
For intuition, the blue marked state nodes in Fig.~\ref{network_ex_step} are the SSC nodes.

\begin{itemize}
\item 
\textbf{Step 1}: Since $\bar{\mathcal{G}}_{1}^{S}=\mathcal{G}_{1}^{S}\cup\mathcal{G}_{12}^{S}$ is \textit{a \textit{tree-type}},
two additional external input nodes $u_1,u_2\in\mathcal{V}_{1}^{IE}$ need to be connected at node $1$ and $3$ to satisfy \textit{Theorem~\ref{theorem_tree}}.
After that, all state nodes of $\bar{\mathcal{G}}_{1}$ are SSC nodes as shown in Fig.~\ref{network_ex_step}(b), 
i.e., $\bar{\mathcal{V}}_{1}^{\text{\tiny SSC}}=\{1,2,3,4,6\}$.\newline

\item 
\textbf{Step 2}: $\bar{\mathcal{G}}_{2}^{S}=\mathcal{G}_{2}^{S}\cup\mathcal{G}_{23}^{S}$ is \textit{a \textit{cycle-type}}.
It follows from \textit{Theorem~\ref{theorem_TCC}} that the component input nodes of $\bar{\mathcal{G}}_{2}$ are $2,3 \in \bar{\mathcal{V}}_{2}^{IC}$.
Hence, one external input node $u_3\in\mathcal{V}_{2}^{IE}$ needs to be connected to node $5$, $7$ or $9$ to satisfy \textit{Theorem~\ref{theorem_cycle}} (we choose node $7$).
After that, all state nodes of $\bar{\mathcal{G}}_{2}$ are SSC nodes as shown in Fig.~\ref{network_ex_step}(c), 
i.e., $\bar{\mathcal{V}}_{2}^{\text{\tiny SSC}}=\{4,5,6,7,8,9,10,12\}$.\newline

\item 
\textbf{Step 3}: Since $\bar{\mathcal{G}}_{3}^{S}=\mathcal{G}_{3}^{S}\cup\mathcal{G}_{34}^{S}$ is \textit{a \textit{cycle-type}},
two input nodes are needed to satisfy \textit{Theorem~\ref{theorem_cycle}}.
According to \textit{Theorem~\ref{theorem_TCC}}, the component input nodes of $\bar{\mathcal{G}}_{3}$ are $4,9\in\bar{\mathcal{V}}_{3}^{IC}$.
Hence, $\bar{\mathcal{G}}_{3}$ already satisfies \textit{Theorem~\ref{theorem_cycle}} without additional external input nodes
since there exists an edge $(10,12)\in\bar{\mathcal{E}}_{3}$ between the nodes $10,12\in\mathcal{N}(\bar{\mathcal{V}}_{3}^{IC})$.
Therefore, all state nodes of $\bar{\mathcal{G}}_{3}$ are SSC nodes as shown in Fig.~\ref{network_ex_step}(d), 
i.e., $\bar{\mathcal{V}}_{3}^{\text{\tiny SSC}}=\{10,11,12,13\}$.\newline

\item 
\textbf{Step 4}: $\bar{\mathcal{G}}_{4}^{S}=\mathcal{G}_{4}^{S}$ is \textit{a \textit{cycle-type}},
two input nodes are needed to satisfy \textit{Theorem~\ref{theorem_cycle}}.
According to \textit{Theorem~\ref{theorem_TCC}}, the component input node of $\bar{\mathcal{G}}_{4}$ is $12\in\bar{\mathcal{V}}_{4}^{IC}$.
Thus, one external input node $u_4\in\mathcal{V}_{4}^{IE}$ needs to be added at node $14$ or $16$ to satisfy \textit{Theorem~\ref{theorem_cycle}} (we choose node $14$).
After that, all state nodes of $\bar{\mathcal{G}}_{4}$ are SSC nodes as shown in Fig.~\ref{network_ex_step}(e),
i.e., $\bar{\mathcal{V}}_{4}^{\text{\tiny SSC}}=\{13,14,15,16\}$.
\end{itemize}

Finally, all the state nodes in $\mathcal{V}^{S}$ are the SSC nodes, i.e., $\bar{\mathcal{V}}^{\text{\tiny SSC}}=\mathcal{V}^{S}$, 
where $\bar{\mathcal{V}}^{\text{\tiny SSC}}={\bigcup}^{4}_{i=1}\bar{\mathcal{V}}_{i}^{\text{\tiny SSC}}$,
From \textit{Corollary~\ref{corollary_TCC}}, it follows that the output graph shown in Fig.~\ref{network_ex_step}(f) is SSC with minimum number of external input nodes.

\begin{rem}\label{remark_complexity}
The \textit{Algorithm~\ref{graph_composition}} provides
the graph-theoretic method to ensure the strong structural controllability of $\mathcal{G}(T)=(\mathcal{V},\mathcal{E})$ with the minimum number of external input nodes.
Similarly, for structural controllability, a method of finding the minimum number of leaders (inputs) has been developed in \cite{guan2021structural}, and
the complexity of Theorem 2 in \cite{guan2021structural} is $\mathcal{O}(|\mathcal{V}|+|\mathcal{E}|)$.
However, the complexity of \textit{Algorithm~\ref{graph_composition}} is $\mathcal{O}(m)$,
where $m$ is the number of disjoint components constituting a pactus. 
Note that $m$ is equal to or less than $|\mathcal{V}|$ because the simplest structure of a pactus is a path, 
e.g., each disjoint component in pactus has only one state node, see \textit{Definition~\ref{definition_sym_pactus}}.
\end{rem}

\section{Conclusion} \label{sec_conc}
This paper delves into the strong structural controllability of undirected graphs in diffusively-coupled networks with self-loops. 
Initially, we establish necessary and sufficient conditions for strong structural controllability, using the concepts of dedicated and sharing nodes. 
Following this, we reinterpret existing results on the controllability of basic components through composition rules applied to disjoint path graphs. 
This approach provides insights into how control paths are topologically configured. 
Using the composition rule, our findings are then extended to pactus graphs, which comprise basic components. 
Finally, we present an algorithm for addressing the MIP for pactus graphs, based on the concept of a component input node, which has the same properties as an external input node.
As our future work, the research presented here can be expanded to include observability, considering its duality with controllability. 
Observability in structured networks is critical, especially in the context of determining the minimal number of measurements necessary for state estimation in diffusively-coupled networks. 
These extensions are earmarked for future exploration.

\section*{Acknowledgment}
This work was supported by both the Institute of Planning and Evaluation for Technology (IPET) of Korea under grant 1545026393 and the National Research Foundation of Korea (NRF) grant funded by the Korea government (MSIT) (2022R1A2B5B03001459).

\bibliographystyle{unsrt}
\bibliography{references}


\balance                              
\end{document}